\documentclass[twoside]{article}
\makeatletter
\let\@fnsymbol\@arabic
\makeatother

\newcommand{\TITLE}{Silent error detection in numerical time-stepping schemes}
\usepackage{amssymb,amsmath}
\usepackage{graphicx}

\usepackage{algorithm}
\usepackage{algpseudocode}
\usepackage{comment}

\newcommand{\bigoh}[1]{\mathcal{O}\left(#1\right)}
\newcommand{\deltat}{{\Delta}t}
\newcommand{\deltax}{{\Delta}x}
\newcommand{\deltay}{{\Delta}y}

\newcommand{\A}{\mathcal{A}}
\newcommand{\B}{\mathcal{B}}

\newcommand{\LapOp}{\nabla \cdot \nabla}



\usepackage[hmarginratio=1:1,top=32mm,columnsep=20pt]{geometry} % Document margins
\usepackage[hang, small,labelfont=bf,up,textfont=it,up]{caption} % Custom captions under/above floats in tables or figures
\usepackage{booktabs} % Horizontal rules in tables
\usepackage{float} % Required for tables and figures in the multi-column environment - they need to be placed in specific locations with the [H] (e.g. \begin{table}[H])
\usepackage{epstopdf}

\usepackage{lettrine} % The lettrine is the first enlarged letter at the beginning of the text
\usepackage{paralist} % Used for the compactitem environment which makes bullet points with less space between them

\usepackage{abstract} % Allows abstract customization
\usepackage{titlesec} % Allows customization of titles
\titleformat{\subsection}
{\normalfont\large\itshape}{\thesubsection}{1em}{}

\title{\vspace{-15mm}\fontsize{24pt}{10pt}\selectfont\textbf{\TITLE}} % Article title

\author{\large
Austin R. Benson\thanks{Institute for Computational and Mathematical Engineering, Stanford University, Stanford, California, USA. (\texttt{arbenson@stanford.edu, schmit@stanford.edu}).}
 \hspace{0.1cm}\footnotemark[2]
\and Sven Schmit\footnotemark[1]
\and Robert Schreiber\thanks{HP Labs, Palo Alto, California, USA (\texttt{rob.schreiber@hp.com}).
We thank the US Department of Energy, which supported this work under Award Number DE - SC0005026.}
\vspace{-5mm}
}
\date{}

\begin{document}

\maketitle % Insert title

%\thispagestyle{fancy} % All pages have headers and footers

%----------------------------------------------------------------------------------------
%	ABSTRACT
%----------------------------------------------------------------------------------------

\begin{abstract}
Errors due to hardware or low level software problems, if detected, can be fixed
by various schemes, such as recomputation from a checkpoint.
{\em Silent errors} are errors in application state that have escaped
low-level error detection.
At extreme scale, where machines can perform astronomically many operations per second,
silent errors threaten the validity of computed results.

We propose a new paradigm for detecting silent errors at the application level.
Our central idea is to frequently compare computed values to
those provided by a cheap checking computation, and to build error detectors based on the difference between the two output sequences.
Numerical analysis provides us with usable checking computations
for the solution of initial-value problems in ODEs and PDEs, arguably the most
common problems in computational science.
Here, we provide, optimize, and test methods based on
Runge-Kutta and linear multistep methods for ODEs, and on implicit and explicit
finite difference schemes for PDEs.
We take the heat equation and Navier-Stokes equations as examples.
In tests with artificially injected errors,
this approach effectively detects almost all
meaningful errors, without significant slowdown.
\end{abstract}

%----------------------------------------------------------------------------------------
%	ARTICLE CONTENTS
%----------------------------------------------------------------------------------------

%\begin{multicols}{2} % Two-column layout throughout the main article text

  \section{Silent errors and checking schemes}
\label{sec:intro}

\begin{comment}
concern for silent error

other ideas in error resilience

fault error failure.
detected and silent errors.
errors that do or do not cause failures.
detected failures in the sense of a troublesome sanity check of results.
detection as a key issue

basic idea.   error can be caused by underresoltion.   A posteriori checks.
For ODEs, the embedded scheme ides, and |A - B| as an approx to the LTE.
Intro the A/B notation.  The whole schema
A/B, D (vector), E (detector), R (response).

For systemic, refine the mesh - traditional.
Question we address: can this idea be used also to detect errors due to the platform.
For PDEs and ODEs?   With implicit methods?
Our goal is to see if we can detect otherwise silent errors at modest cost.
R out of scope here.  A - D in scope.

response.  recomputed solves transient
change of method/resolution solves systemic error
change of machine solves hard errors

Detecting errors: problems of false positives and false negatives.
Adjusting the sensitivity to control these.

Into schemes:
A/B:  RK explicit
Adams explicit, Implicit.
Full LMMs including extrapolations
PDE - require implicit, eg Crank-Nicholson (essential an Adams Moulton method).

Intro diff:  Just the infinity norm for now.

Detectors:   Use diff relative to nearby diffs as a check.  Self adjusting.   No need to specify a threshold.
Into the two detector approach.   Into the self adjustment of threshold.   Need experiments to show the effects of these.

Look at the benefit of an unstable A scheme.

\end{comment}

\subsection{Silent errors are worrisome}
Computational scientists are concerned about silent errors in exascale computing.
Silent errors are perturbations to application state
that may lead to a failure such as a bad final solution \cite{Snir:2013}.
These errors may arise from a bit flip, a firmware bug, data races, and other causes.
Several authors (\cite{Cappello:2009, Dongarra:2011, Snir:2013}) have discussed
the sources and the frequency of silent errors.

Why the current concern?
An exaflop machine will be able to do on the order of $10^{23}$ operations per day, and will have on the order of
$10^{17}$ bytes of memory~\cite{Dongarra:2011}.  And in order to achieve very aggressive energy efficiency and performance targets, machine architects are pushing envelopes:
with near threshold voltage logic, with new memory and storage technologies, and with photonic communication.
Consumer quality hardware may already suffer errors at the personal computer scale once per year~\cite{Nightingale2011},
and cost precludes really significant hardening of the hardware in supercomputers.
Thus, the scale of systems makes such errors quite likely.
Indeed, some high-performance systems today already suffer from silent errors at a troublesome rate~\cite{shi2009testing}.

\subsection{Algorithmic responses to silent errors}
The numerical algorithms community has already looked at error vulnerability.
It is well known that many errors do not cause failures.
Other errors lead to an obvious application failure.
Silent errors are more worrisome, because they can
cause unsuspected erroneous outputs.   Our goal is to make these errors non-silent.

It has been argued that with extra care, convergent iterations are inherently self-correcting; for example,
a resilient version of GMRES is proposed in \cite{Hoemmen:2011}.
Other empirical studies have shown, however, that iterative methods are sometimes vulnerable to errors \cite{Bronevetsky:2008, Casas:2012}.
And in a study of a minimization approach to Hartree-Fock ground state calculation,
Van Dam, {\em et al.}, found that ``it is
insufficient to rely on the algorithmic properties of the
Hartree-Fock method to correct all the possible bit-flips and
resulting data corruption''~\cite{vanDam:2013}.
Sufficiently big errors can be fatal to these algorithms.

Minimization and equation solving (in which the data defining the function or equation are assumed to be incorruptible) is an easy case, since the residual of the current approximate solution almost surely does not lie.  But in many cases in computational science, a time-dependent, initial value problem is solved.  In these cases, any perturbation to the computed solution puts that solution onto a permanently erroneous track.
We therefore take the view that error detection is a fundamental issue, and that errors, detected as soon as they occur,
can then be handled by an appropriate correction scheme.

Certain common kernels have been fortified with error detectors.
For example, checksum methods have been used for matrix multiplication \cite{Huang:1984},
high performance LU \cite{Du:2012}, and checking the integrity of
data replicated on multiple compute nodes \cite{vanDam:2013}.
This latter paper also monitored a number of theoretical invariants of the algorithm.
Orthogonality of a matrix can be checked by multiplying by its transpose, for example.
Conservation laws can be monitored, where available, as a check for error.
These monitoring approach were found to be useful, but fallible; they are not a comprehensive safety net.

\subsection{Our approach}
In this paper, we propose a very general, low-cost error detection approach
that applies to iterative computations in general, and to the solution of initial value problems for ODEs and PDEs in particular.
Our central idea is to compare the solution given by a primary, or {\em base} time stepping scheme to the solution given by an auxiliary {\em checking} scheme, and to do this every time step.  The two schemes use the same input data, all of it computed by the base scheme: thus the auxiliary solution is used only locally, at each time step, to check for errors.
This approach has compelling advantages over a straightforward duplication of the computation -- it is cheaper, and it can detect problems that duplication cannot.

Error is a constant in scientific computing.
Even with no bugs or failures,
we have modeling error, truncation error, and roundoff error.
%the former the purview of applied math and engineering, the latter two the concern of numerical analysis.
To deal with
truncation error, schemes of the kind we are proposing have long been employed for automatic step size control in
ODE solvers~\cite{Fehlberg:1969}.  Similar {\em a posteriori} error estimators are used for mesh adaptation in PDE solvers~\cite{Berger1983}.
The idea is to make a numerical method introspective, aware of and watchful for errors.
Our contribution is to extend these powerful schemes so that
they can be used to detect errors due to a misbehaving computing system as well.

Suitable checking schemes are available in the most common
setting -- the solution of initial-value problems in ODEs and PDEs.
A full description of the general approach, including the key question of how we trigger notification of an error,
is given in the next section.
We describe specific checking schemes for Runge-Kutta and linear multistep methods in Section~\ref{sec:ode_solvers} and
for finite difference methods for the heat equation in Section~\ref{sec:PDE_solvers}.
We discuss the error detector, and an approach to controlling and optimizing it, in Section IV.

In tests with artificially injected errors,
we measure the impact of errors by how much they impact the solution.
We quantify this idea in Section~\ref{sec:inject},
and we then show  through numerical experiments that our detection scheme effectively catches errors that have a significant impact on the solution.

  \section{Outline of general method}
\label{sec:method}

Suppose there are two iterative methods to solve a problem, a base method, $\B$, and an
auxiliary checking method, $\A$.
One would use $\B$ in a computing environment with no errors.
Desirable properties of $\B$ are therefore accuracy and stability.
A suitable auxiliary $\A$ solves the same problem -- its output can be compared to that of $\B$.
And it can be used at each iteration, using the same input data as does $\B$. % S: this is correct right?
The key idea is that the norm of the difference between the results provided by $\A$ and $\B$ is an estimator of the magnitude of the difference at the current step between $\B$ and an error-free solution.
That suggests that $\A$ should have accuracy comparable to (or even better than) $\B$.
For efficiency, we want $\A$ to be fast when used as a check on $\B$, possibly by reusing some of the computations, communications, and
input data of $\B$.
Since we do not use $\A$ in a closed-loop setting (i.e., $\A$ does not use its own results as input at each step),
stability is not really important in $\A$.  This gives us useful freedom in choosing auxiliary schemes.

The two schemes produce sequences of values $\{A_i\}$ and $\{B_i\}$ in the same normed vector space.  For any norm or seminorm of the difference, we have a scalar
sequence, $D_i = \| A_i - B_i \|$.
We may choose to use more than one such metric, so in general, $D_i$ may be a vector.

Our methods employ, at each step $n$, a window into the sequence
$(D_{n - d}, \ldots, D_n)$, as data for an error detection function $E(D_{n - d}, \ldots, D_n)$
that decides whether or not to raise the flag for an error.
The error detector typically employs one or more measures of size, or more powerfully measures of anomaly, to the
value $D_n$ in the context of its recent values $D_{n-d}, \ldots, D_{n-1}$.
In general, then, a method includes base and checking computation schemes, a vector of difference measures, and an
error detection criterion.

\subsection{Choosing an auxiliary scheme}
Let us first consider the simplest case, and show why it does not help us;
this motivates the search for better auxiliary schemes.

Consider $\A$ to be the same as $\B$, i.e., $\A$ just repeats the computation.
The error detector simply flags an error whenever $\B$ and $\A$ differ, or differ by more than a small multiple of machine precision.
Here $\A$ is not fast (it costs the same as $\B$).
Moreover, it won't catch certain errors: if the input data for the step are corrupted (after successful previous steps
write these data to memory), the two computations produce identical, incorrect results.
The same is true if a computation unit fails in a repeatable manner -- a stuck-at fault, for example -- or a communicated value corrupted by the network is used by both schemes.

A better checking scheme is one that reuses some of the computation and communication of the base scheme (for efficiency),
but is different in a way that makes the two schemes disagree unless there are no errors anywhere in the algorithm.

Two examples, covered in detail below, are embedded Runge-Kutta schemes,
in which $\A$ reuses evaluations of the derivative function that are needed for $\B$,
and paired linear multistep methods, in which saved and new values of the solution
and its time derivative are combined in two different ways to estimate the solution at
the next time step.
Notably, for stiff ODEs and second-order parabolic PDEs, stability mandates the use of implicit base schemes,
with their attendant algebraic system to solve at each time step,
but explicit schemes prove to be effective, inexpensive auxiliaries.

\subsection{The detector}
The norm of the difference, $D_i = \| A_i - B_i \|$, is an obvious candidate for error checking.
What complicates this is that the size of $D_i$ can vary over orders of magnitude in the error-free case, as the solution changes.
Thus, a hard threshold is ineffective for error detection.
Our view is that a sudden change in the sequence $\{D_i\}$ better indicates that an error is present.
We examine this empirically in Section~\ref{sec:numerical}.

At the core of the detector there will be comparisons of some scalar indicator quantities to some thresholds.
How should these thresholds be chosen?  We take the following view.
With any sort of error detector, we can have false positives and false negatives,
and there is an intrinsic tradeoff between their rates:
a low threshold boosts the rate of false positives but misses few true
errors; a high threshold reduces the false positive rate at the expense of more missed errors.
In our setting, false negatives are bad, but false positives are only mildly annoying.
We will propose below strategies that tend to
maintain a constant and sufficiently small false positive rate, while reducing the false negative rate to very low levels.

\begin{comment}
Our method for finding errors depends on two assumptions:
\begin{enumerate}
  \item \emph{Methods $\B$ and $\A$ should give approximately the same result in the error-free case.}
  \item \emph{The difference between their results should normally change only slowly.}
\end{enumerate}
The first means that the schemes are checks of one another.   The second allows us to appropriately choose a threshold
for the difference, so as to minimize missed errors (false negatives) without an excessive rate of false positives. Section~\ref{sec:examples} illustrates these issues.
%This second observation does not always hold;
%the step size should be sufficiently small and it also
%relies on the continuity of the solution.
%However, it gives reason to believe that the approximations of $\B$ and $\A$ should
%remain of similar size in consecutive iterations.
% Hence, if we combine these two ideas, they can be used to predict whether a silent error has occurred.
\end{comment}

\subsection{What to do if an error is flagged}
\label{sec:whatiferrorflagged}
What do we expect an application to do if an error is detected?
This is not the topic of our work, but we feel it's important to give
an idea of a general scheme.

Before implementing any error recovery scheme for an iterative method, there are two important questions to ask when a flag is raised to indicate an error:
\begin{enumerate}
    \item On what iteration could the error have occurred?
    \item What data or computation was affected by the fault?
\end{enumerate}
We then intend to redo the failed steps (Question 1) by redoing all potentially failed computations (Question 2).

How far back must we go?
In our numerical experiments, we see that in quite a few cases the time step \emph{after} the one in which the fault occurred causes the error flag to be raised.
For example, a small error in a derivative evaluation in a linear multistep method may produce a large error in a few time steps, since
the derivative evaluation gets re-used for several iterations.
Thus, we first have to establish which iterations may be erroneous and which ones we still trust.

For example, suppose the base and checking schemes use
the solutions at the previous two time steps to compute the solution at the next.
Futher suppose that this pair of schemes sometimes flags an error one iteration after the faulty step. Consider now that our procedure flags iteration $5$, signalling there might be an error.
We cannot trust the solution at step $4$, but we can trust step $3$ because, were it faulty,
we would have seen a flag at iteration $3$ or $4$.
Hence, we go back and restart the computation of iteration $4$, using the stored data from iterations $2$ and $3$.

We hope and expect that silent errors will be rare.
Even with a small false positive rate, there will be
more false positives than true positives.
It is important to not get stuck and redo (correct)
computations over and over again when they are incorrectly flagged.
In order to avoid this we propose
a taxonomy of possibilities.

\begin{enumerate}
\item
The error may have been caused by a transient fault.
On retry, if the fault does not recur, we will likely succeed, with no error flag.

\item
The error may have been caused by a permanent fault that causes erratic, irreproducible behavior.

\item
First, if data in memory have been corrupted, silently, we may discover this with our scheme.
If we redo the failed steps starting from the same corrupted in-memory data,
we expect an identical outcome, with the error flag raised on the retry.

\item
Second, if some hardware or software component has failed in a ``hard'' way,
meaning it consistently produces incorrect results, we expect again an identical outcome,
with the error flag raised on the retry.

\item
The error may be expected if, for example, the algorithm sacrifices correctness in rare scenarios for speed \cite{Rinard:2013}.
If the algorithm fails randomly, but rarely (as is the case in \cite{Rinard:2013}), re-computation will likely provide the correct answer.
If the failure is deterministic, then a different algorithm will be needed.

\item
Or there may be no error at all; it may be that the problem's local difficulty causes the error flag to trigger,
with the current value of the time step, in the error free case.
Again, we expect an identical result on retry.

\end{enumerate}

As stated above, the first response to a flagged error is
to go back to an iteration that was computed with no error flags, or possibly one iteration further back, and redo the subsequent iterations, including the error checks.
We want also to check to see whether this recomputation produces the same result as it did initially.

If retry succeeds with no error flag, we likely have discovered that an error of type (1) or (5) occurred.
We can report it, and continue.

If retry fails, but with a significantly different result that on the first try, it is likely that a component of the platform has become unreliable; we need to change it.
This is an error of type (2).

If retry fails with the same computed result as initially, there is ambiguity: the cause may be any of (3), (4), (5), or (6) above.
Absent a way to tell, there is a problem.

Our approach is compatible with systems that protect memory contents from loss and corruption.
All machines have a basic error detection and correction system for memory.
These can be augmented with additional protection and redundancy, as in the Global View Resilience Project~\cite{GVR}.
In case of a reproducible flagged error, such a scheme can be invoked to test the memory contents and see if they have been corrupted, to rule our an error of type (3),
or correct it if one has occurred.

If such a test detects no corruption of the input data to the failed step, then how can we distinguish between errors of type (4) and (6)?
Here we could go back to the origins of this approach, and redo the computation with a smaller step -- perhaps half the current step.
If the error is of type (6), this approach will, after a few reductions, fix the problem.
But if the problem remains after repeated step size reductions, we would suspect an error of type (4).

  \section{Applications}
\label{sec:examples}

\subsection{ODE solvers} \label{sec:ode_solvers}

\begin{figure*}
  \centering
  \fbox{\includegraphics[scale=0.7]{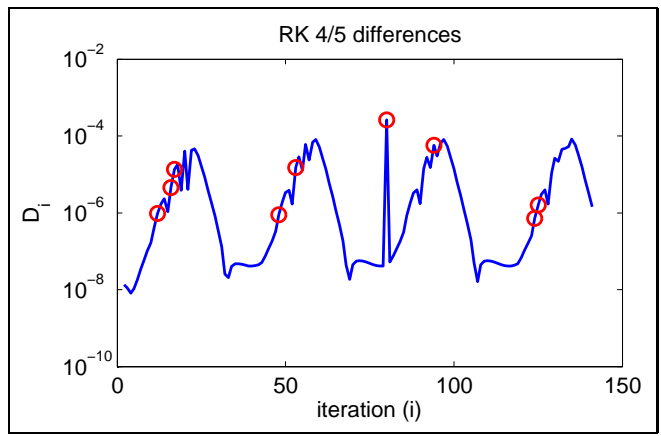}}
  \caption{
    Difference between RK5 ($\B$) and RK4 ($\A$) over time for the van der Pol equation with
    $b = 2$ and initial conditions $u(0) = 1$ and $u^{\prime}(0) = 0$.
    An artificial error is injected at the $80$-th iteration, which results in the spike in $D_{80}$.
    The red circles indicate iterations that are predicted to be erroneous by our detection scheme; see Section~\ref{sec:error}.
  }
  \label{fig:rk45}
\end{figure*}

Consider a first-order ODE initial-value problem:
\begin{equation}
  \frac{d}{dt}u(t) = f(t, u(t)), \quad u(0) = u_0.
  \label{eq:ode}
\end{equation}
%% Runge Kutta

Suppose that we are using the explicit midpoint Runge-Kutta (RK) scheme as a base scheme $\B$ to compute $u^{\B}_{n+1} \approx u(t_{n+1})$ in Equation~(\ref{eq:ode}):
\begin{subequations}
  \begin{align}
     k^{\B}_1 &= f(t_n, u^{\B}_n)  \label{eq:RKMP1} \\
     u^{\B}_{n+1} &= u^{\B}_n + hf(t_n + \frac{1}{2}h, u^{\B}_n + \frac{1}{2}hk^{\B}_1) \label{eq:RKMP2}.
  \end{align}
\end{subequations}

The local truncation error (LTE) of this scheme is $\bigoh{h^3}$.
We note that $f(t_n, u^{\B}_n)$ is the central computation in Euler's method, which has LTE $\bigoh{h^2}$, and we use this to construct an auxiliary scheme $\A$,
\begin{equation}
u^{\A}_{n+1} = u^{\B}_n + hk_1^{\B}. \nonumber
\end{equation}

An example difference computation is $D_{n+1} = \| u^{\B}_{n+1} - u^{\A}_{n+1} \|_{\infty}$.
By re-using $u_N^{\B}$ and $k_1^{\B}$, $\A$ provides a cheap approximation to the solution.
The midpoint and Euler schemes are an {\em embedded RK pair} \cite{Dormand:1980}; these form the basis of adaptive step-size methods.
In general, we can use any embedded RK pair in the $\A/\B$ formulation.
A common, accurate scheme is the RKF45 scheme due to Fehlberg~\cite{Fehlberg:1969}.

Figure~\ref{fig:rk45} illustrates how errors lead to jumps in the difference between $\A$ and $\B$ using this particular scheme for the Van der Pol equation
\begin{equation} \label{eq:vdp}
u^{\prime\prime}(t) - b(1 - u(t)^2)u^{\prime}(t) + u(t) = 0,
\end{equation}
whose rapid changes in derivatives make this a challenging case.

In Section~\ref{sec:numerical}, we show in more detail that errors in the evaluation of Equation~(\ref{eq:RKMP1}) or Equation~(\ref{eq:RKMP2}) can
be effectively detected by RK-based $\A/\B$ schemes.
Moveover, we also show that errors in the evaluation of $f$ can be detected just as effectively,
\emph{even if this wrong computation is used by both methods}.
However, due to the memoryless property of RK methods, they cannot easily detect changes to $u_n^{\B}$.
We discuss this matter and provide experiments in Section~\ref{sec:vdp}.

%% Multistep
Linear multistep methods (LMM) are also amenable to our framework.
An Adams-Bashforth LMM (AB-LMM) of order $p \geq 1$
computes $u^{\B}_{n+1} \approx u(t_{n+1})$ by
\begin{equation}
u^{\B}_{n+1} = u^{\B}_n + \sum_{i=n-p+1}^{n}{h\alpha_{p,i}f(t_i, u^{\B}_i)}, \nonumber
\end{equation}
such that the LTE is $\bigoh{h^{p+1}}$.
Suppose that $\B$ is a $p$-th order AB-LMM, $p \ge 2$.
One choice of $\A$ is the AB-LMM of order $p-1$,
which reuses the same data, stores no additional data, and performs no additional evaluation of $f$.
An alternative $\A$ is a linear multistep method of order $p$ that interpolates at (possibly multiple) $u^{\B}_{k}$ for $k < n$.
However, additional memory is needed to store solutions at prior time steps.
In order to compare AB-LMM to RK methods, we will consider the ($p-1$, $p$) AB-LMM pairs in our experiments in Section~\ref{sec:numerical}.

Implicit numerical schemes are preferred for stiff ODEs.
For example, an Adams-Moulton LMM (AM-LMM) of order $p$
defines $u_{n+1}$ implicitly, as the solution to the system of equations
\begin{equation}
u^{\B}_{n+1} = u^{\B}_n + \sum_{i=n-p+2}^{n+1}{h\beta_{p,i}f(t_i, u^{\B}_i)}. \nonumber
\end{equation}
Its LTE is $\bigoh{h^{p+1}}$.
Suppose that the base scheme $\B$ is an AM-LMM.
The computationally expensive part of the method is solving the (possibly nonlinear) equation for $u_{n+1}$.
A lower-order AM-LMM will require a different solve and is a less attractive choice for $\A$ (it will not be fast compared to $\B$).
Instead, we can use an AB-LMM for $\A$,
\begin{equation}
u^{\A}_{n+1} = u^{\B}_n + \sum_{i=n-p+1}^{n}{h\alpha_{p,i}f(t_i, u^{\B}_i)} \nonumber
\end{equation}
AB-LMM is an explicit method, but the starting value, $u^{\B}_n$, and function evaluations, $\{f(t_i, u^{\B}_i)\}_{n-p+1 \le i \le n}$, have been computed by the implicit AM-LMM.
Thus, use of an AB-LMM as $\A$ does not suffer from the instability of explicit methods used on stiff ODEs.
Note that we can employ any implicit LMM as a base scheme, including, for example, the backward differentiation formulas.

Finally, an explicit / implicit pair of LMMs is sometimes used in a predictor-corrector fashion.
Here, the implicit LMM's equations are solved by a truncated fixed-point iteration in which the explicit scheme generates a first iterate.
In this instance, the auxiliary scheme (the predictor) is already a part of the solution mechanism for the base scheme (the corrector), so it comes at no extra cost.

\subsection{PDE solvers}\label{sec:PDE_solvers}

\begin{figure*}[tbp]
  \centering
  \fbox{\includegraphics[scale=0.7]{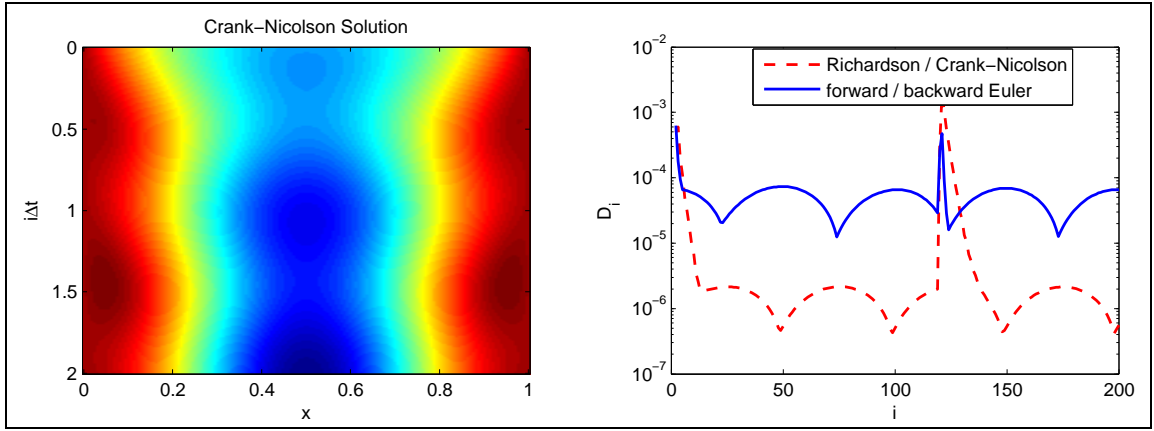}}
  \caption{
    Solution and difference sequence for Equation~(\ref{eq:heat}) with
    $k = \frac{1}{100}$, $q(x, t) = 0.1\left(\sin(2\pi t) + \cos(2\pi x)\right)$, $v(x) = x(x-1)$,
    $\deltax = \frac{1}{160}$, and $\deltat = \frac{1}{100}$.
    The difference function is $D_i = \| A_i - B_i \|_{\infty}$.
    A fault is injected at the 120-th time step by multiplying the 40-th component of the right-hand-side of the linear solves used by Crank-Nicolson and backward Euler by $0.995$.
    The solution appears normal, but the difference sequence indicates an anomoly at the 120-th time step.
  }
  \label{fig:heat_soln_diffs}
\end{figure*}

Due to the large variety of PDE solvers, we do not have a one-size-fits-all solution.
For time-dependent PDEs, a method of lines discretization in space results in a system of ODEs,
and the ODE methods described above can be employed.
Here instead we consider finite difference schemes for PDEs.

To make this idea concrete, we will describe an $\A/\B$ formulation for the heat equation.
A more detailed example (the incompressible Navier-Stokes equations) is provided in Section~\ref{sec:navierstokes}.

For a model problem, consider the nonhomogeneous heat equation
\begin{align}
& u_t = ku_{xx} + q(x, t), \quad k > 0 \nonumber \\
& u(x, 0) = v(x)\label{eq:heat}
\end{align}
with homogeneous Dirichlet boundary conditions.
Suppose that $\B$ and $\A$ are the backward and forward Euler schemes.
Both methods have LTEs of $\bigoh{(\deltax)^2}$ in space and $\bigoh{\deltat}$ in time.
At each time step, $\B$ solves a linear system, while
$\A$ computes a matrix-vector product.
Thus, we expect $\A$ to be faster.
Moreover, in the distributed memory setting, $\A$ requires no communication other than what is done in scheme $\B$, if $\B$ uses an iterative solver.

An alternative $\B$ is the Crank-Nicolson scheme, which is implicit and has LTEs of $\bigoh{(\deltax)^2}$ in space and $\bigoh{(\deltat)^2}$ in time.
Forward Euler is a candidate for $\A$, but we desire an explicit method with the same LTEs as Crank-Nicolson.
The Richardson scheme (also known as the leapfrog scheme) is such a method, and it uses a centered difference in time and space.
While the Richardson scheme reuses the computations from Crank-Nicolson, the centered difference in time requires the solution at the two previous time steps.
We refer to \cite[Section~6.3]{Strikwerda:2007} for a discussion of all of these methods.

Figure~\ref{fig:heat_soln_diffs} plots the solution to the heat equation and the sequence of differences $\{ D_i = \| A_i - B_i\|_{\infty} \}$ for a particular problem instance.
An error occurs in the 120-th iteration and is exhibited by a spike in the sequence of differences.
In Section~\ref{sec:numerical_heat}, we thoroughly examine how effective the Richardson / Crank-Nicolson and forward / backward Euler $\A/\B$ formulations are in detecting errors.

\subsection{Extrapolation}
For some base schemes $\B$, the choice of a related auxiliary $\A$ may not be obvious.
But extrapolation is always available, in the form of
an LMM in which all the $\beta$ terms are zero.
The order 1 version is simply
\begin{equation}
  A_i =  2B_{i-1} - B_{i-2}. \nonumber
\end{equation}
Although useless as basic solvers, extrapolation methods are usable for error detection.
They are cheap, and can have acceptable error characteristics.
Extrapolation is part of the $\A/\B$ scheme used for the Navier-Stokes equations in Section~\ref{sec:navierstokes}.

  \section{Error detection}
% Might want to remove "and resolution"
\label{sec:error}

In this section we outline a practical implementation of the error detection function, $E$.
Throughout this section, we consider a scalar error metric $D$, for example the sup-norm of the difference
between $A$ and $B$.
Because the scale of the variations in $D$ is unknown, and can be time varying,
it is important that $E$ can handle errors independent of scale.
We therefore use relatively large changes in $D$ as indicators of problems, rather
than the magnitude of $D$ itself.

To find iterations with errors,
we use two indicator variables derived from the sequence of differences $D_{n-p}, \ldots, D_{n}$:
\begin{align}
  J_n &= \frac{D_{n+1} - D_{n}}{D_{n}} \label{eq:relative_jump} \\
  V_n &= \frac{\text{Var}(D_{n-p+1},\ldots D_{n+1})}{\text{Var}(D_{n-p},\ldots,D_{n})} \label{eq:var_change}
\end{align}
$J_n$ measures the jump in the sequence and $V_n$ measures a change in variance.
A large value for either indicator signals that an error has occurred.
The integer $p$ adjusts the window size; in our experiments, $p = 10$.

We flag an error only when both indicators exceed their current thresholds.
We show in Section~\ref{sec:detect} that the two-indicator strategy improves the sensitivity
to actual errors for a fixed false positive rate.
\begin{comment}
Then, we move back $x$ steps, where $x$ depends on how many iterations might be erroneous,
as discussed in Section~\ref{sec:whatiferrorflagged}.
\end{comment}

We use a closed-loop mechanism to tune the thresholds.
We increase a threshold by the factor $\Gamma > 1$
every time the indicator is above its threshold
and decrease a threshold by a factor $\gamma < 1$
every time the indicator is below its threshold.
If both indicators are above their respective thresholds, we flag an error and rely on an error handler as discussed in Section~\ref{sec:whatiferrorflagged}. 
%to determine
%whether or not the error is a false positive.
The overall idea is in Algorithm~\ref{alg:resilient}.

The closed-loop tuning procedure forces thresholds to be only slightly above typical values.
In practice, this reduces the probability of false negatives, without causing many false positives.
In our experiments in Section~\ref{sec:numerical}, $\Gamma = 1.4$, $\gamma = 0.95$, and $p=10$.
The detector's performance is not sensitive to these choices.
We choose $\gamma$ to be close to one so as to reduce the false positive rate.

%% Add algorithm pseudocode
\begin{algorithm}
\begin{algorithmic}
    \Procedure{Resilient algorithm}{}
      \State Initialize thresholds $\tau_J$ and $\tau_V$
      \State Initialize increase parameter $\Gamma > 1$.
      \State Initialize decrease parameter $\gamma < 1$.
      \While{$n < N$}
        \State $B_{n+1} =$ BaseMethod()
        \State $A_{n+1} =$ AuxiliaryMethod()

        \State $D_{n+1} = \| B_{n+1} - A_{n+1} \|$
        \State \texttt{// Compute indicators}
        \State $J_n = \frac{D_{n+1} - D_{n}}{D_{n}}$
        \State $V_n = \frac{\text{Var}(D_{n-p+1},\ldots D_{n+1})}{\text{Var}(D_{n-p},\ldots,D_{n})}$

        \State \texttt{// Check for errors}
        \If{$J_n > \tau_J$ and $V_n > \tau_V$}
          \State FlagError()
          \State Move backward: $n = n - x$
        \Else
        \State \textsc{UpdateThreshold}($J_n$, $\tau_J$)
        \State \textsc{UpdateThreshold}($V_n$, $\tau_V$)
        \State Move forward: $n = n + 1$
        \EndIf
      \EndWhile
    \EndProcedure
    \State
% Updating thresholds procedure to make the algorithm less cluttered.
    \Procedure{UpdateThreshold}{$t$, $\tau$}
%    \State Consider $\Gamma>1$ and $\gamma < 1$ fixed.
      \If{$t > \tau$}
        \State $\tau = \Gamma \tau$
      \Else
        \State $\tau = \gamma \tau$
      \EndIf
    \EndProcedure
  \end{algorithmic}
  \caption{
  Pseudocode for the error detection algorithm.
  We use adaptive error thresholding to keep the false positive rate and false negative rate low.
  }
  \label{alg:resilient}
\end{algorithm}

\subsection{Discussion}
Our method relies on two assumptions:
\begin{itemize}
  \item Methods $\B$ and $\A$ produce approximately the same result in the error-free case; that is, they are accurate.
  \item Changes in the solution are not excessively rapid, so that the behavior of the difference sequence is predictable in the error-free case.
\end{itemize}
When either of these assumptions is violated, we expect problems.
For example, a discontinuity in the derivative $f$ in Equation~(\ref{eq:ode}) may cause consecutive iterations to vary wildly, and an error will often be predicted.
And indeed we may have detected an error: not one caused by an unstable platform, but rather one due to underresolution.
It is thus inherent in our approach that when the solution changes rapidly, we may
see false positives.

We note that sometimes an error is flagged one step \emph{after} the iteration the error occurred.
It is therefore best to go back two iterations.
We explore this phenomenon in Section~\ref{sec:tardy}.

Many more sophisticated statistical tools are available in time series analysis for outlier and peak detection~\cite{Hamilton:1994, Lin:2003}.
However, algorithms for peak detection are not a focus of this paper, and
the simple and fast approach outlined above performs well for several practical examples in Section~\ref{sec:numerical}.

  \section{Numerical experiments}
\label{sec:numerical}

\begin{comment}
In this section we discuss the numerical results obtained from the resilient algorithm.
Faults are injected like this
No timing results due to implementation in Matlab -> figure mathmatically.
For some results regarding the heat equation, see Figures~\ref{fig:heat2}~and~\ref{fig:heat3}.
We notice, as expected, that smaller timesteps imply better detection rates.
Also, we note that small errors are much harder to detect.
\end{comment}

% make this a little longer
In this section, we evaluate the performance of several $\A/\B$ formulations on a variety of problems.
We begin by outlining how we inject artificial faults in computations.
Next, we elaborate on the particular problem instances and how the detector performs.
All experiments used Matlab R2013a.

\subsection{Fault injection and LTE-normalized error}
\label{sec:inject}

In our experiments, we inject faults by corrupting important computations or data,
such as the result of a function evaluation.
We do not corrupt internal data structures or program logic.
The reason for this choice is that low-level errors in the program will likely either \emph{not} be silent or will
have a similar effect to corrupting computations or data.
The three ways we inject faults are as follows:
\begin{enumerate}
\item Corrupt an evaluation of $f$, the derivative function in Equation~(\ref{eq:ode})
 or an evaluation of $q$, the source term in Equation~(\ref{eq:heat}).
\item Corrupt the right-hand-side when solving a system of linear equations \emph{before} the solver is used.
\item Corrupt a previous solution from the solver (most of our solvers use the solution from the previous time step).
\end{enumerate}

By ``corrupt'', we mean multiply a single component of a vector or matrix by some amount.
In our experiments, we conduct many trials and multiply by a normally distributed random variable with mean $1$ and problem-dependent variance, $\sigma^2$.
This method avoids producing \texttt{Inf}s and \texttt{NaN}s that are produced by fault-injection methods such as bit flips.
Such values are easily recognized and flagged, so the errors are not silent.

Suppose that a fault is injected at iteration $n-1$.
In order to measure the impact of a fault on the computed solution
relative to the ordinary truncation errors of the numerical method, we use the value
\begin{equation}
L_n = \frac{\| B_n - \hat{B}_n \|}{\| \hat{B}_n - \hat{A}_n \|}, \nonumber
\end{equation}
where $\hat{B}_n$ and $\hat{A}_n$ are the outputs of $\B$ and $\A$ when no fault is injected.
The numerator measures the magnitude of the impact of the fault on the base solution, and the denominator is an estimate of the LTE.
We call $L_n$ the \emph{LTE-normalized error}.

By measuring the LTE-normalized error, the type of error introduced in the experiments becomes less important.
Instead, we are able to see how detection varies with the error's impact on the solution.
A small LTE-normalized error means that the error has relatively little influence on the solution while a large LTE-normalized error means that the error has a large influence on the solution.
In the subsequent sections, we show that our detector effectively catches large LTE-normalized errors but has difficulty catching small LTE-normalized errors.
In practice, this is a desirable property: errors that have a stronger impact on the solution are more easily detected.
Furthermore, it seems unreasonable to demand that errors of the order of local truncation error are detected.

\subsection{Van der Pol equation}\label{sec:vdp}
% Dutch names: Capitalize first letter of first and last word, hence
% Robin van Persie or Van Persie
%
% Glad we have a Dutchman on this project--lots of Dutch numerical analysts. -Austin
% And now you will know this forever, and be annoyed every time you see this mistake. -Sven

Our first set of experiments uses the Van der Pol equation (Equation~(\ref{eq:vdp})).
We will vary the damping parameter $b$ in our experiments but fix the intial conditions and time interval:
\begin{equation}
    u(0) = 1, \quad u'(0) = 0, \quad t \in [0, T] = [0, 14] \nonumber
\end{equation}
The van der Pol equation has rapid changes in derivatives, which makes it a difficult test problem for our error detection scheme.
Increasing $b$ stiffens the problem and induces more rapid changes in derivatives.

We test four $\A/\B$ schemes: Runge-Kutta 4/5 (RK45), Runge-Kutta 2/3 (RK23), Adams-Bashforth 4/5 (AB45), and Adams-Bashforth 2/3 (AB23).
In the first set of experiments, we corrupt one component of the derivative evaluations at some time step ($\sigma^2 = $ 1e-1).
Runge-Kutta uses four (RK23) or six (RK45) derivative evaluations per step, and the error corrupts one of these evaluations.
Adams-Bashforth uses one derivative evaluation per step, so the time step determines the corrupted function evaluation.
The erroneous time step, corrupted derivative component, and erroneous evaluation in Runge-Kutta are chosen uniformly at random.
We use 2,000 trials, where a trial consists of an ODE solve at times $0, h, 2h, \ldots, T$.
One error is introduced per trial.
Finally, we use two values of the damping parameter, $b \in \{2, 3\}$.
For $b = 2$, the step sizes are $1 / 10$ and $1 / 20$ for the Runge-Kutta and Adams-Bashforth methods, respectively.
For $b = 3$, the step sizes are $1 / 15$ and $1 / 35$.

\begin{figure*}
  \centering
  \framebox{
  \parbox{400\unitlength}{
  \includegraphics[scale=0.75]{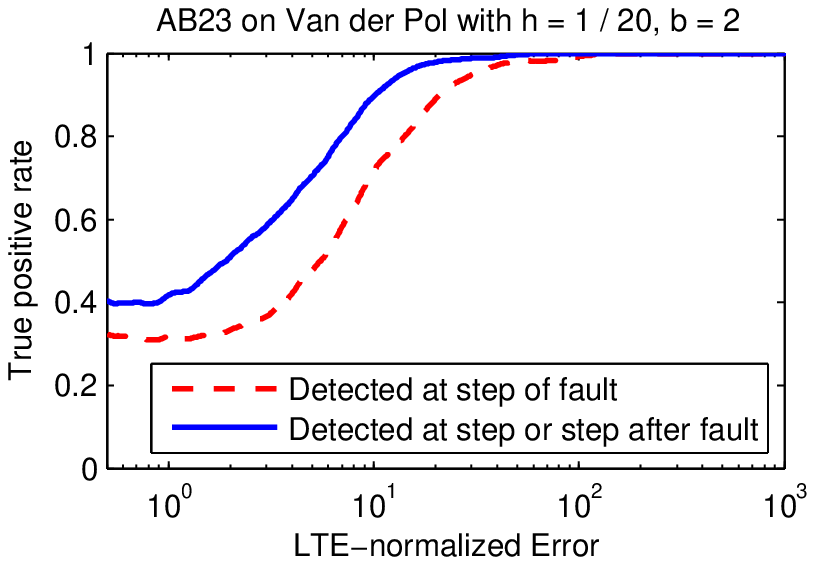}
  \includegraphics[scale=0.75]{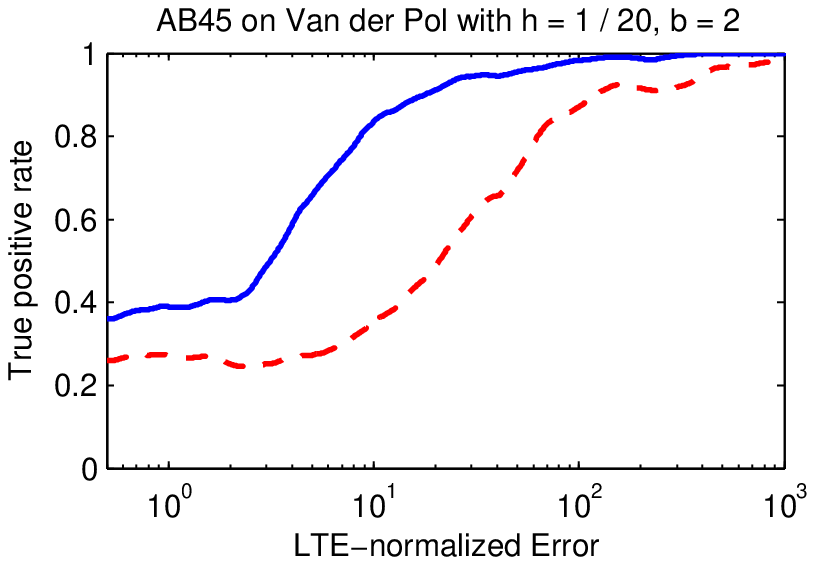} \\
  \includegraphics[scale=0.75]{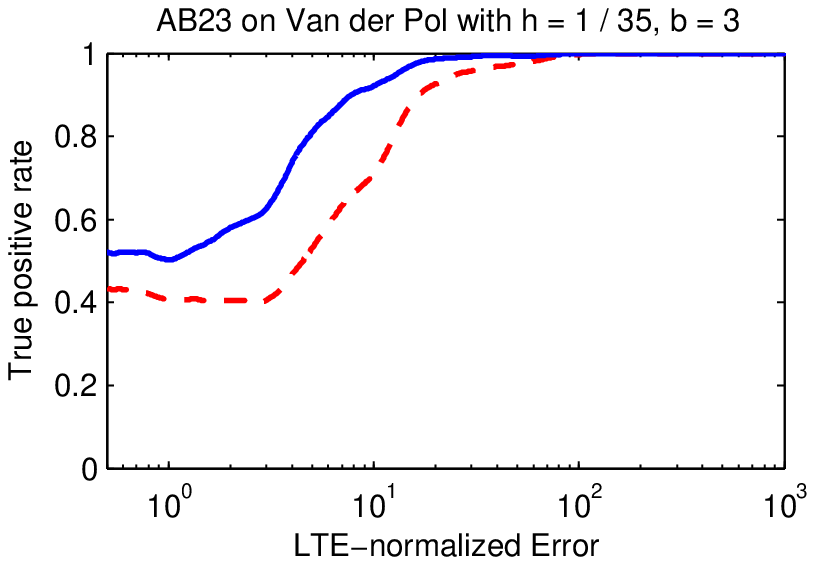}
  \includegraphics[scale=0.75]{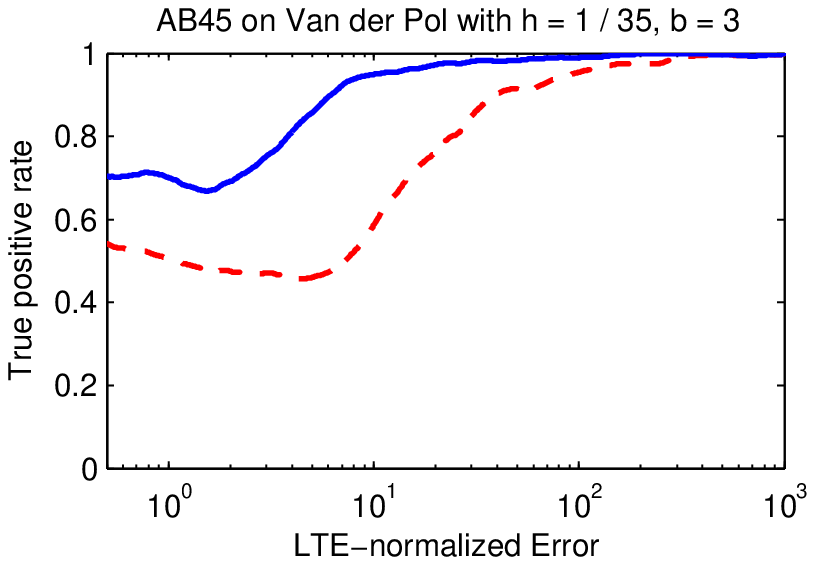} \\
  \includegraphics[scale=0.75]{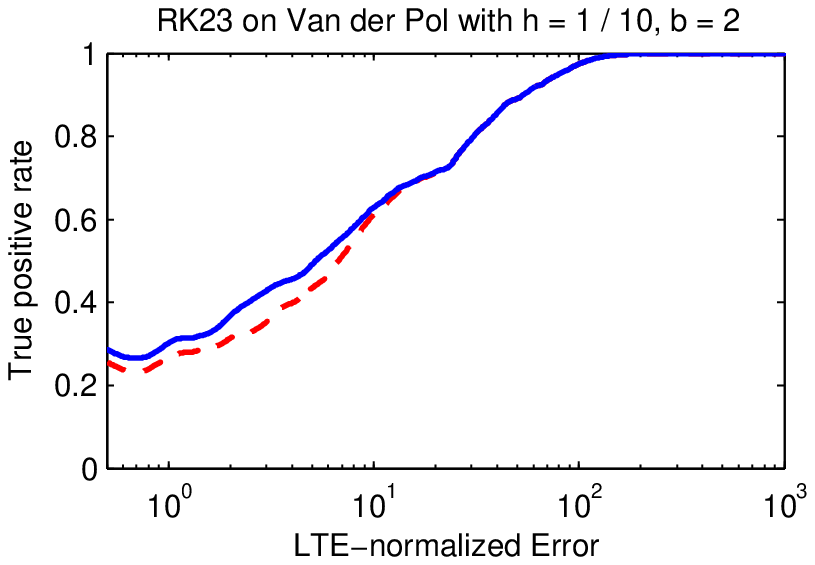}
  \includegraphics[scale=0.75]{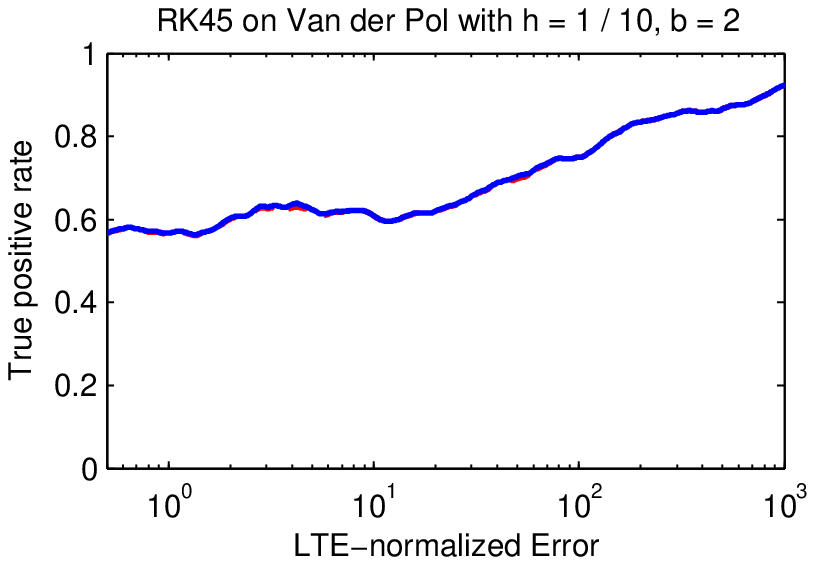} \\
  \includegraphics[scale=0.75]{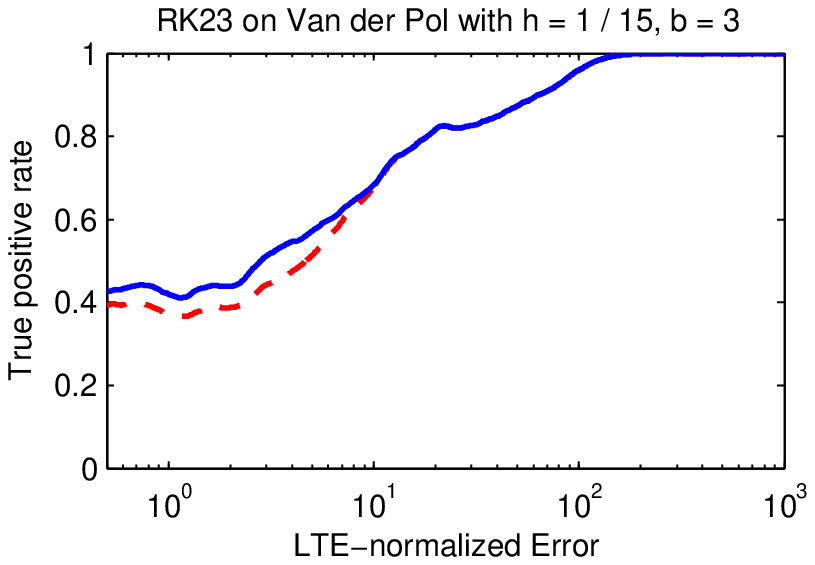}
  \includegraphics[scale=0.75]{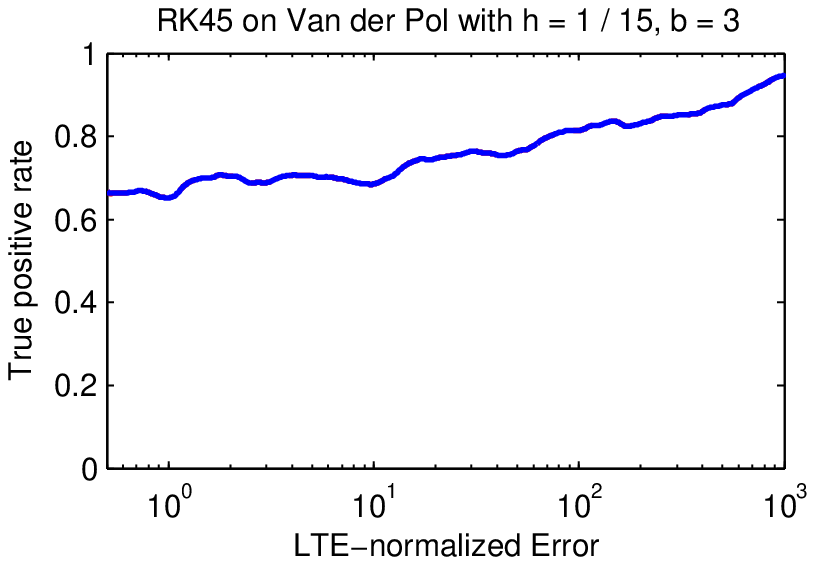}
  }}
  \caption{
    Detector performance with RK45, RK23, AB45, and AB23 $\A/\B$ schemes on the Van der Pol equation.
    We corrupt a single derivative evaluation at the current time step by multiplying one component of the evaluation
    by a normal random variable with mean 1 and variance 1e-1.
    The same sequence of corruption amounts (values of the random variable) was used for each plot.
    Kernel regression with a Gaussian kernel was used to compute the curves.
  }
  \label{fig:vdp_func}
\end{figure*}

Figure~\ref{fig:vdp_func} shows the true positive rate (TPR) as a function of the LTE-normalized error.
The true positive rate is the proportion of artificially injected errors detected by the detection scheme.
We use a kernel regression with a Gaussian kernel to fit the TPR to the LTE-normalized error.
Each plot shows the detection rate for both (1) detection at the time step of the fault and (2)
detection at the step of the fault or the step after.
In all cases, we see the trend that large LTE-normalized errors are easily detected while small LTE-normalized errors are more difficult to catch.
Contrary to RK45 and RK23, AB45 and AB23 detect many errors the step after the error occurs.
This is not entirely surprising.
Runge-Kutta methods use the erroneous derivative evaluation once to advance a time step,
while Adams-Bashforth methods reuses the erroneous computation at the time step of the fault \emph{and} in following steps.
Thus, there is more opportunity for the $\B$ and $\A$ schemes to disagree in Adams-Bashforth methods.
Finally, we note that, in general the higher-order schemes (RK45 and AB45) exhibit slightly better performance than the lower order schemes (RK23 and AB23).

In the second set of experiments, we corrupt a previous time step data stored in memory.
Error-correcting codes in memory hardware provides a low-level check for faults that corrupt data in memory, and applications can supplement these with other, perhaps stronger protections, at some cost; but it is interesting to find whether our approach can detect changes in stored data independently.

For Runge-Kutta, the relevant stored data is the solution computed at the last time step.
For an LMM, the stored state is a set of several solutions and derivatives at previous time steps.

\begin{figure*}
  \centering
  \framebox{
  \parbox{400\unitlength}{
  \includegraphics[scale=0.75]{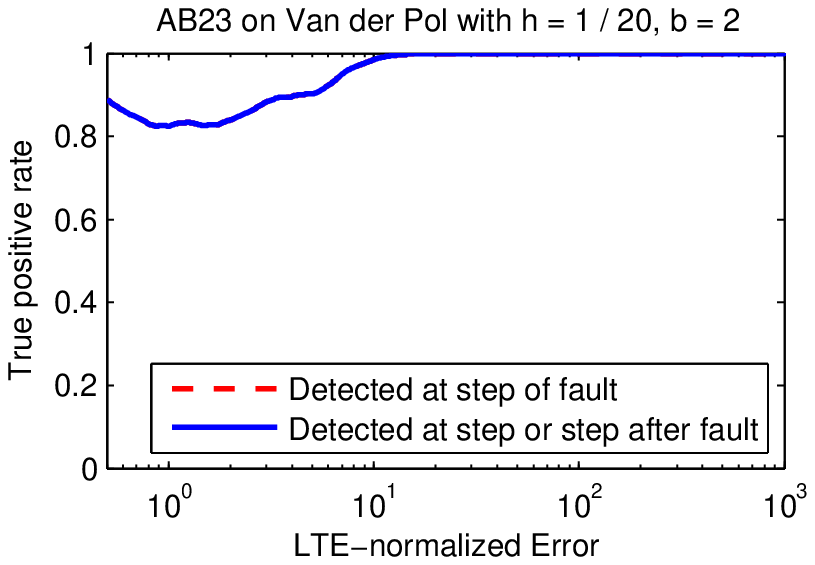}
  \includegraphics[scale=0.75]{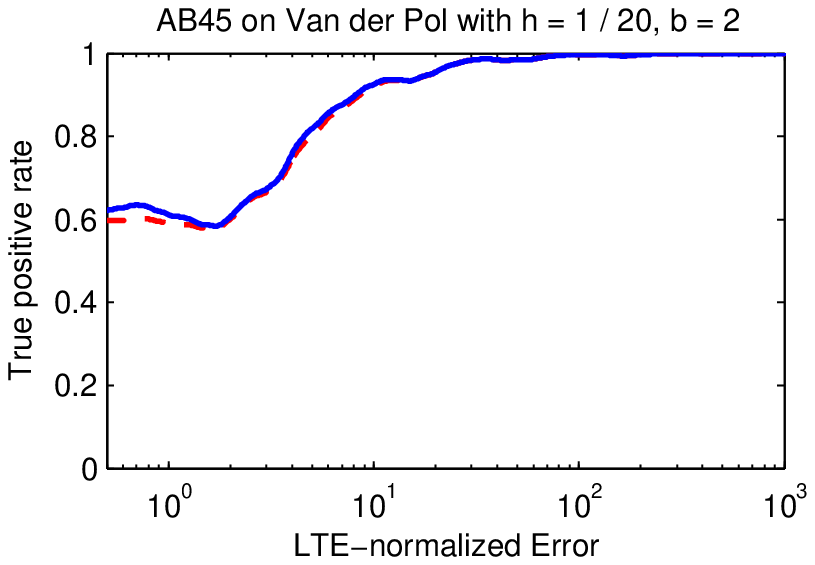} \\
  \includegraphics[scale=0.75]{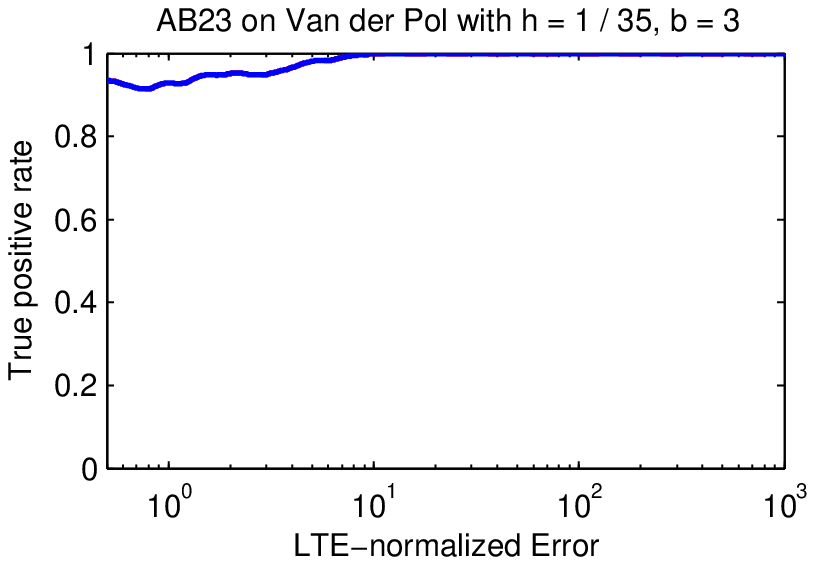}
  \includegraphics[scale=0.75]{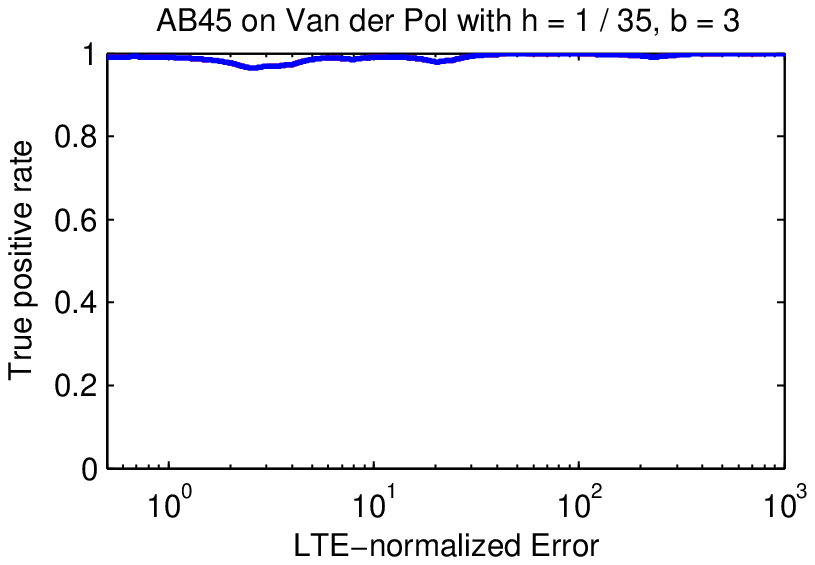} \\
  \includegraphics[scale=0.75]{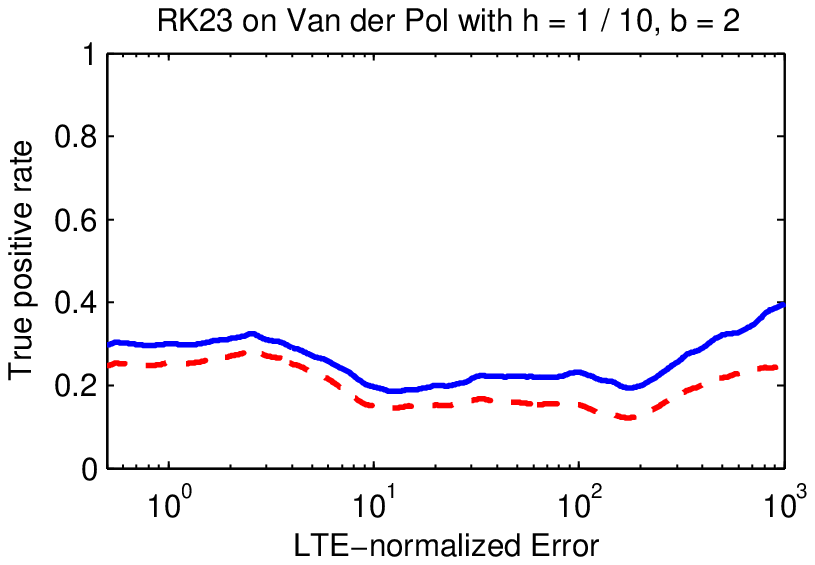}
  \includegraphics[scale=0.75]{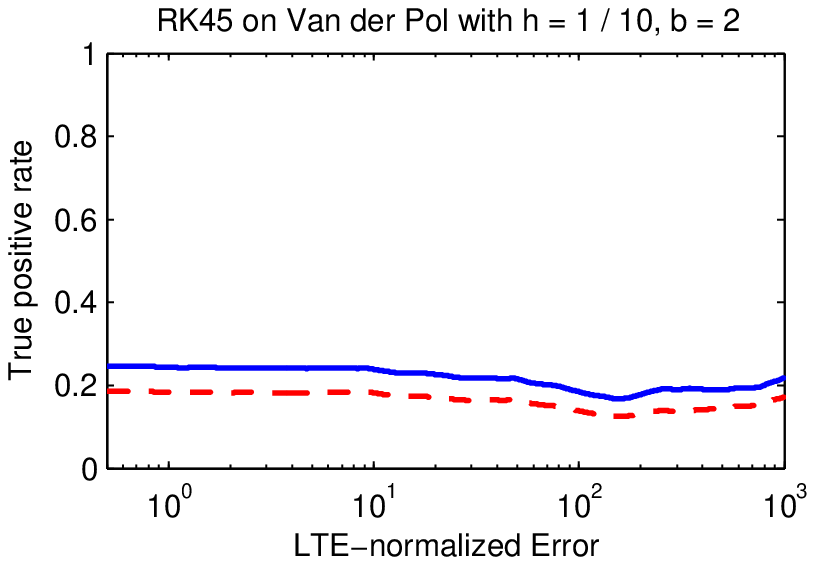} \\
  \includegraphics[scale=0.75]{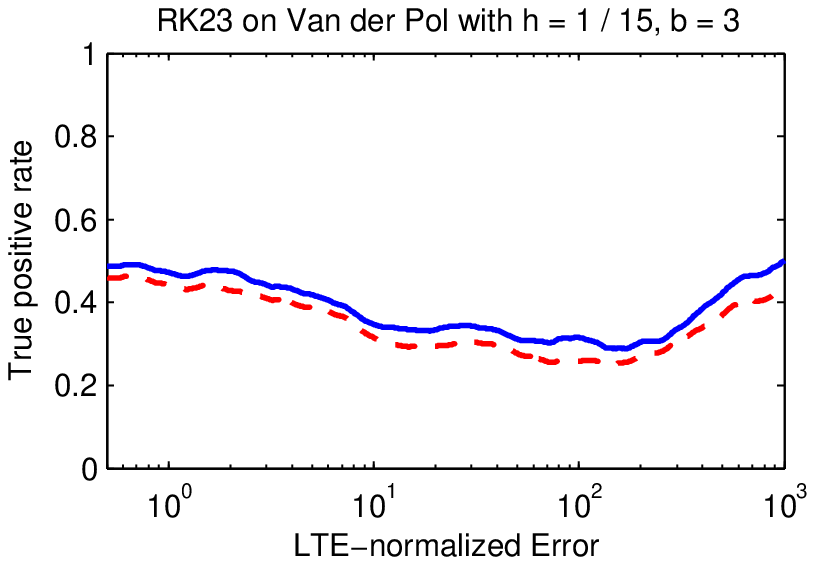}
  \includegraphics[scale=0.75]{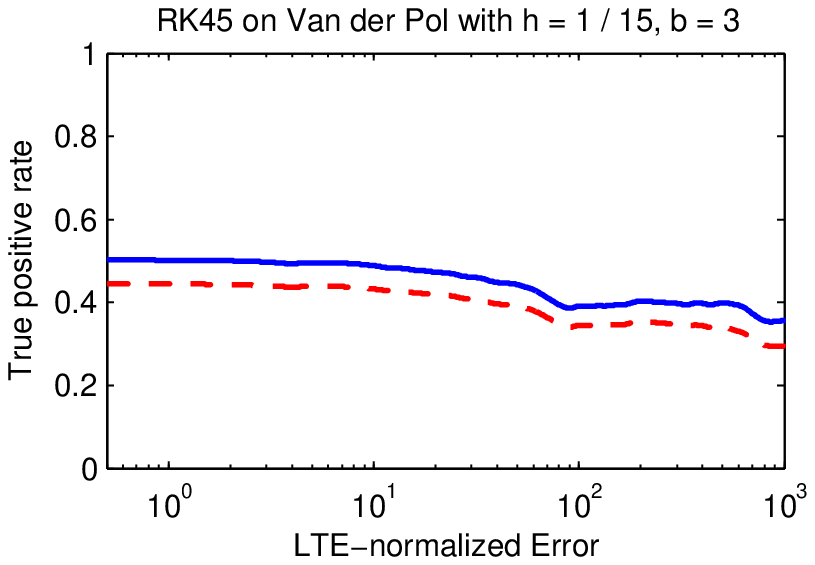}
  }}
  \caption{
    Detector performance with RK45, RK23, AB45, and AB23 $\A/\B$ schemes on the Van der Pol equation.
    We corrupt the solution from the last time step (or a previous derivative evaluation in AB23 and AB34)
    by multiplying one component of the vector by a normal random variable with mean 1 and variance 1e-1.
    The corrupted data is stored in memory in the program.
    The same sequence of corruption amounts (values of the random variable) were used for each plot.
    Kernel regression with a Gaussian kernel was used to compute the curves.
  }
  \label{fig:vdp_prev}
\end{figure*}

Figure~\ref{fig:vdp_prev} shows the error detection effectiveness of Runge-Kutta and LMM-based schemes.
We see that the Runge-Kutta $\A/\B$ schemes have difficulty detecting the errors.
At each step of any Runge-Kutta, the previous solution is the initial condition for advancing to the next time step.
Thus, the difference computation $D_{n} = \| u_n^{\B} - u_n^{\A} \|$ does not necessarily seem out of the ordinary;
$D_n$ is the correct difference for the wrong problem.
The change in initial conditions can cause $D_n$ to be significantly larger than $D_{n-1}$, so the detection rates are still modest.
With a multistep method, on the other hand, the previous solution and derivative evaluations need to be correct for $D_n$ to be the correct difference.
Thus, AB45 and AB23 detect these errors effectively;
the TPR is quite high when the error is the result of corrupting stored solution or derivative data.

These results illustrate an advantage of linear multistep methods compared to one-step methods.
It could be argued that a checksum could be used to detect changes to data stored in memory, and
that these could be used in conjunction with one-step $\A/\B$ schemes for error detection -- one can perform a check on the memory content at each step, before accepting the solution at the next step.
But these schemes cannot detect data corruption due to a bug that stores an incorrect value.
Thus, it is important to be able to detect memory data corruption at the application program, and multistep schemes appear to do this effectively.

\subsection{Heat Equation}\label{sec:numerical_heat}

\begin{table*}[tbp]
\centering
\begin{tabular}{c c c c c c c}
\toprule
Configuration & $q(x, t)$ & $v(x)$ & $k$ & $T$ & $\deltax$ & $\deltat$ \\ \midrule
1 & $xe^{-t/2}$ & $4x(x-1)(x-2)$ & $\frac{1}{100}$ & 2 & $\frac{1}{100}$ & $\frac{1}{60}$, $\frac{1}{100}$, $\frac{1}{140}$ \\
2 & $\frac{1 - \sqrt{1 - 4(t - t^2)}}{2 - 2t}$ & $6 \vert x - \frac{1}{2} \vert - 3$ & $\frac{1}{1000}$ & 1 & $\frac{1}{200}$ & $\frac{1}{100}$, $\frac{1}{200}$, $\frac{1}{400}$ \\
3 & $0.1\left(\sin(2\pi t) + \cos(2\pi x)\right)$ & $x(x-1)$ & $\frac{1}{100}$ & 2 & $\frac{1}{160}$ & $\frac{1}{100}$, $\frac{1}{160}$, $\frac{1}{200}$ \\ \bottomrule
\end{tabular}
\caption{Three configurations of the heat equation.}
\label{tab:heat_configs}
\end{table*}

\begin{figure*}[tbp]
  \centering
  \framebox{
  \parbox{470\unitlength}{
  \includegraphics[scale=0.7]{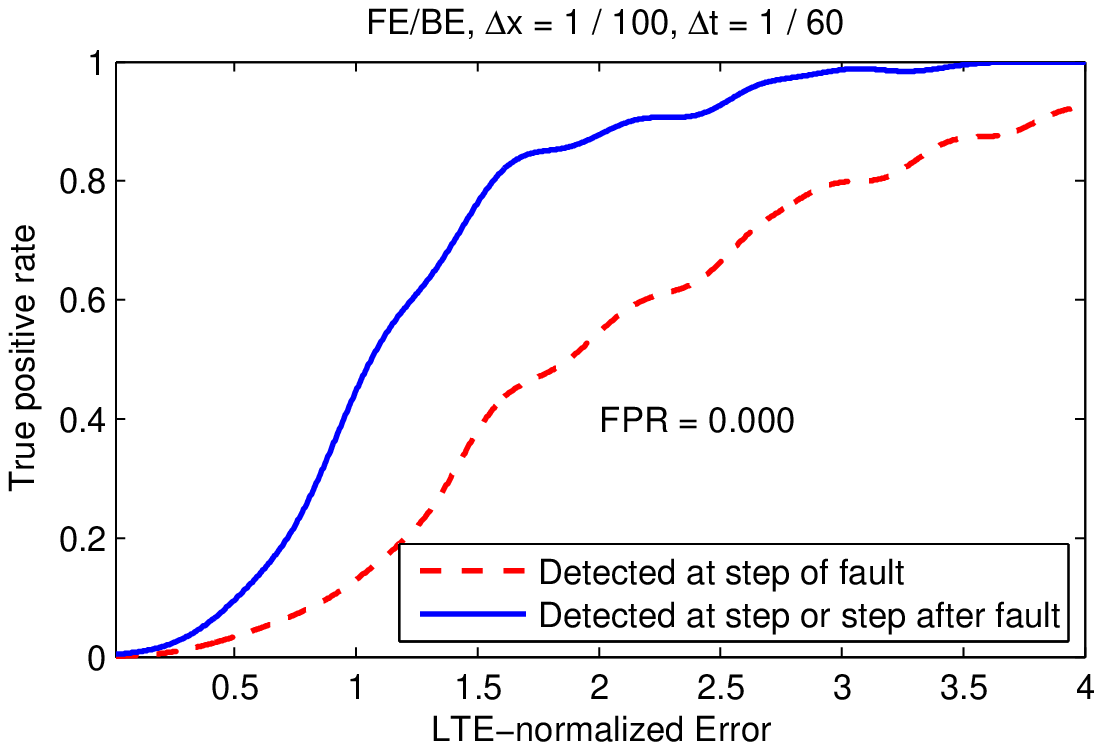}
  \includegraphics[scale=0.7]{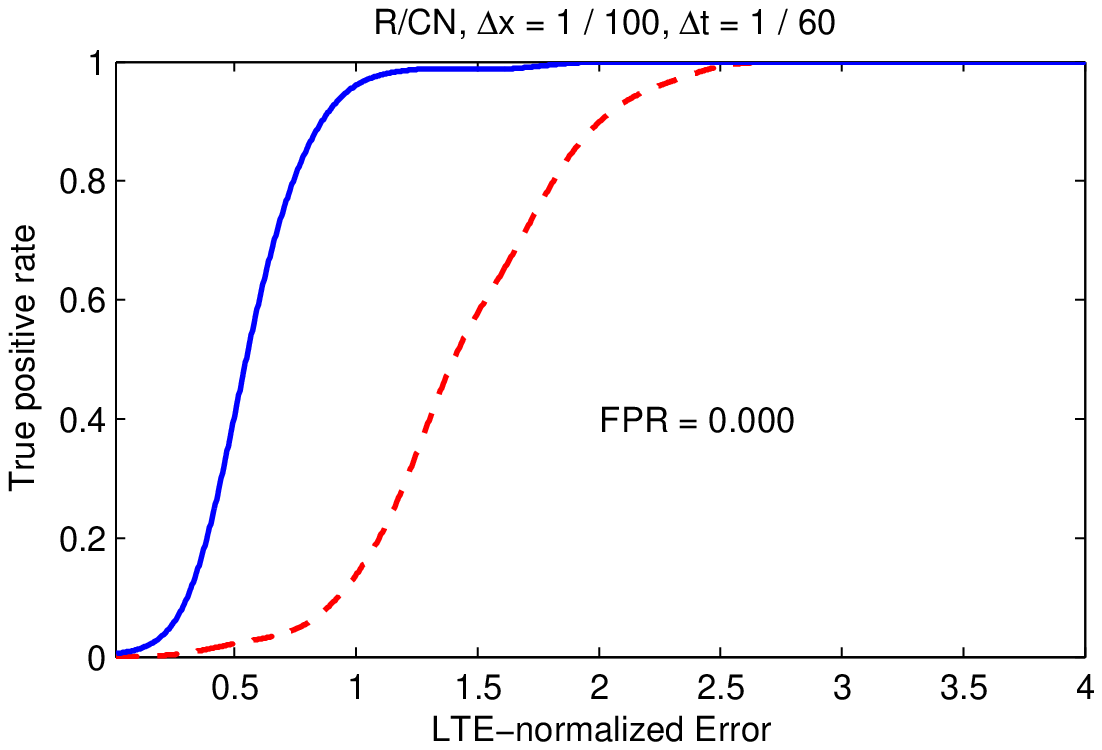} \\
  \includegraphics[scale=0.7]{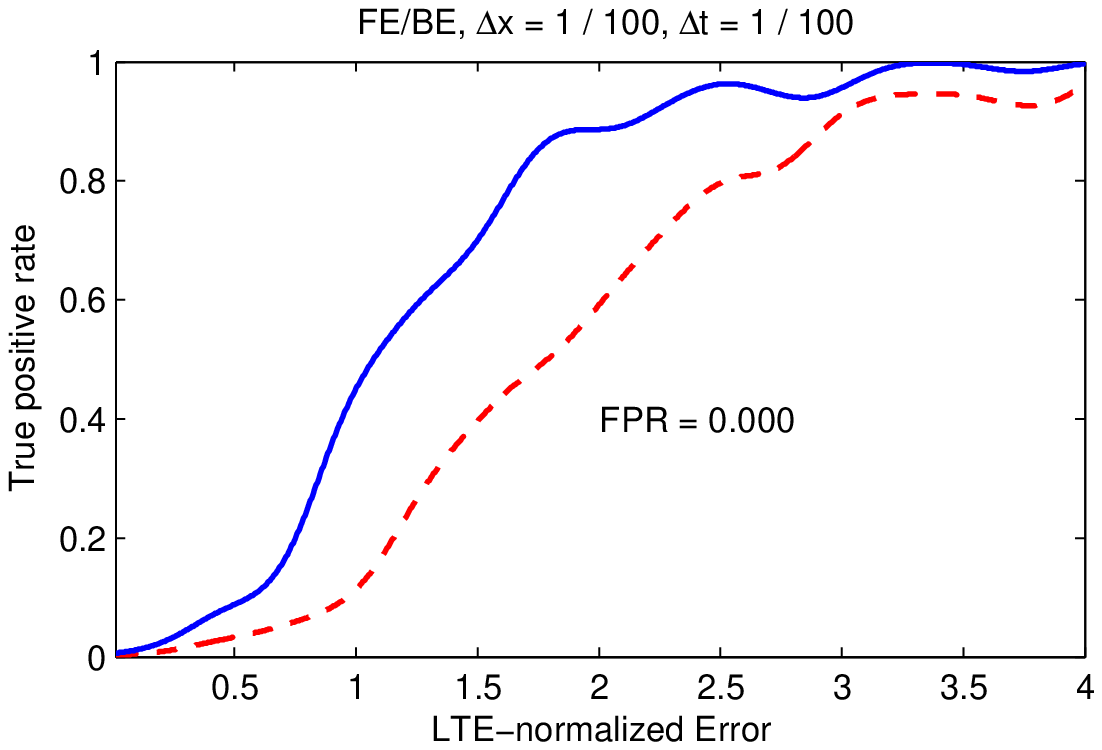}
  \includegraphics[scale=0.7]{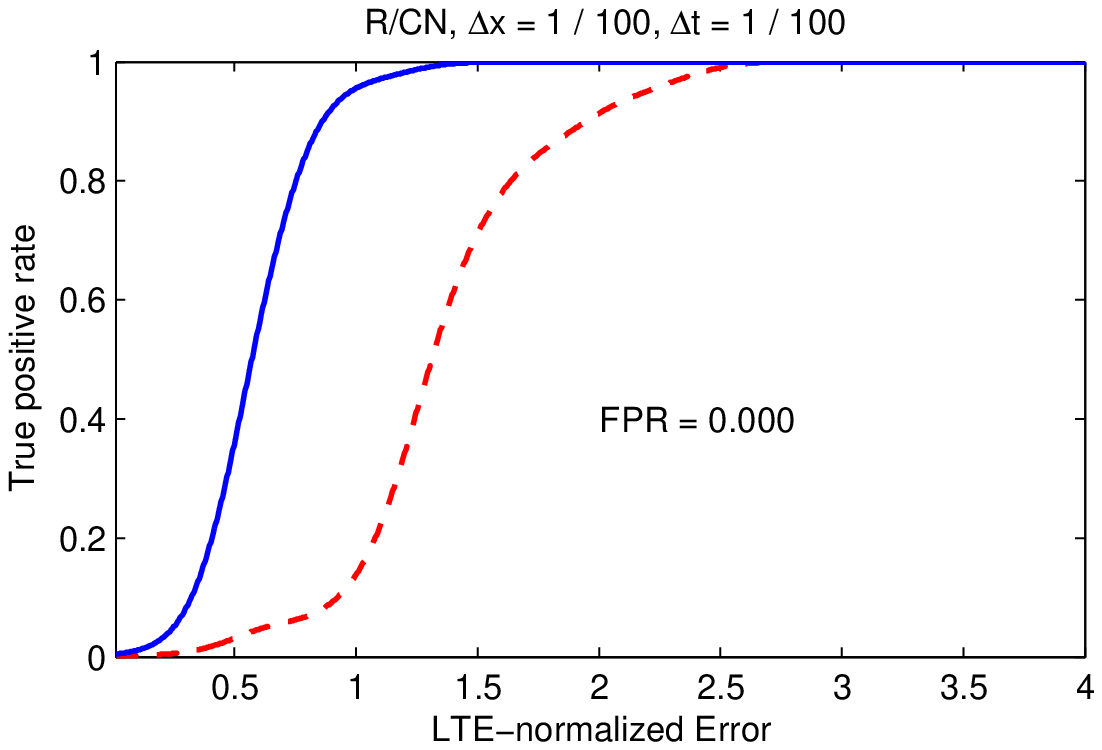} \\
  \includegraphics[scale=0.7]{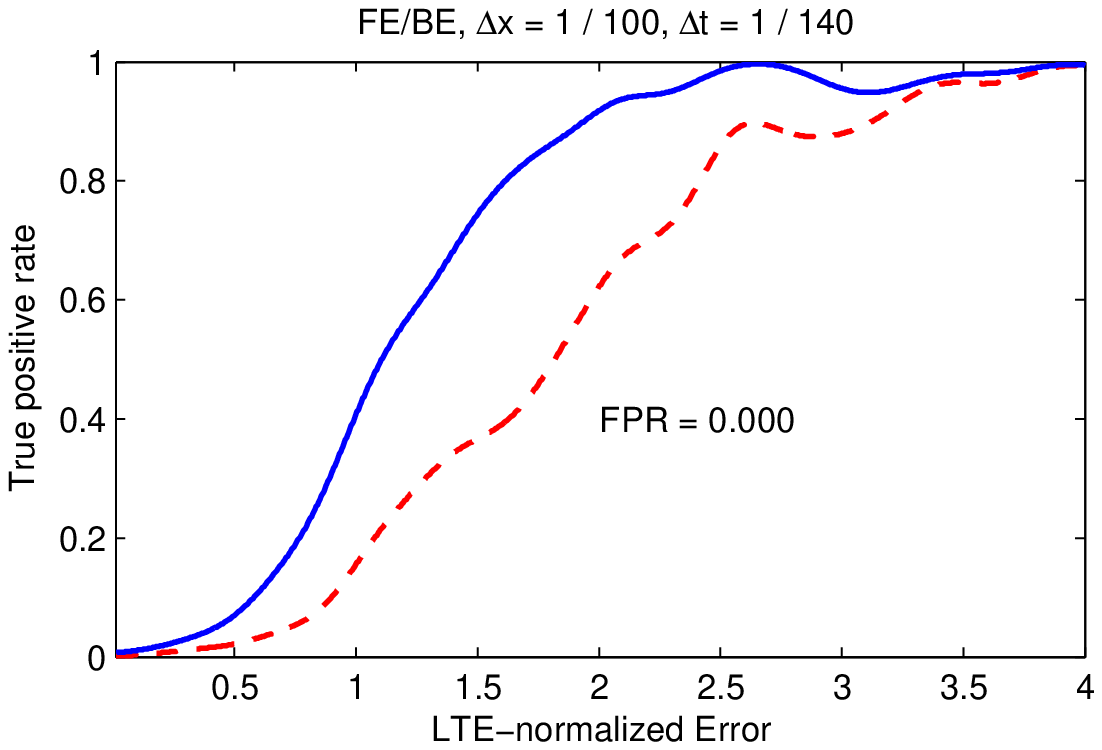}
  \includegraphics[scale=0.7]{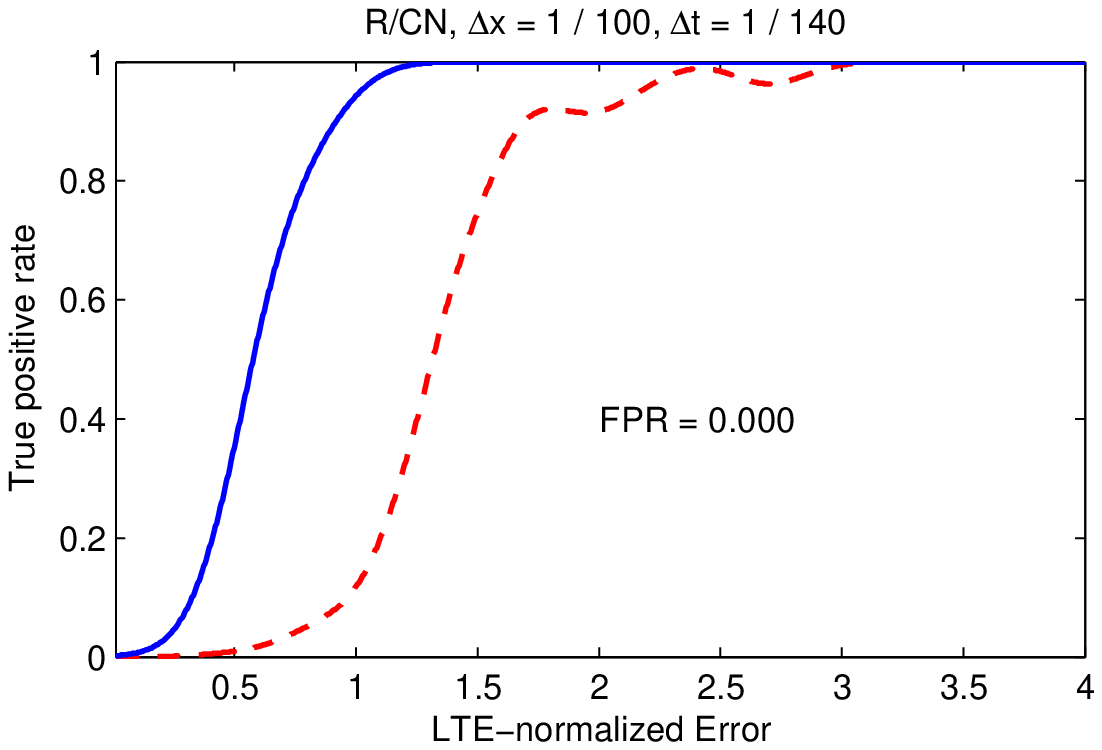}
  }}
  \caption{
  Detection results for Heat equation under Configuration 1 with $dt = 1/60$, $1 / 100$, and $1 / 140$.
  Faults are injected by multiplying the source term $q$ by a normally distributed random variable with
  mean $1$ and variance 1e-3 (R/CN) or 1e-1 (FE/BE).
  }
  \label{fig:heat1}
\end{figure*}

\begin{figure*}[tbp]
  \centering
  \framebox{
  \parbox{470\unitlength}{
  \includegraphics[scale=0.7]{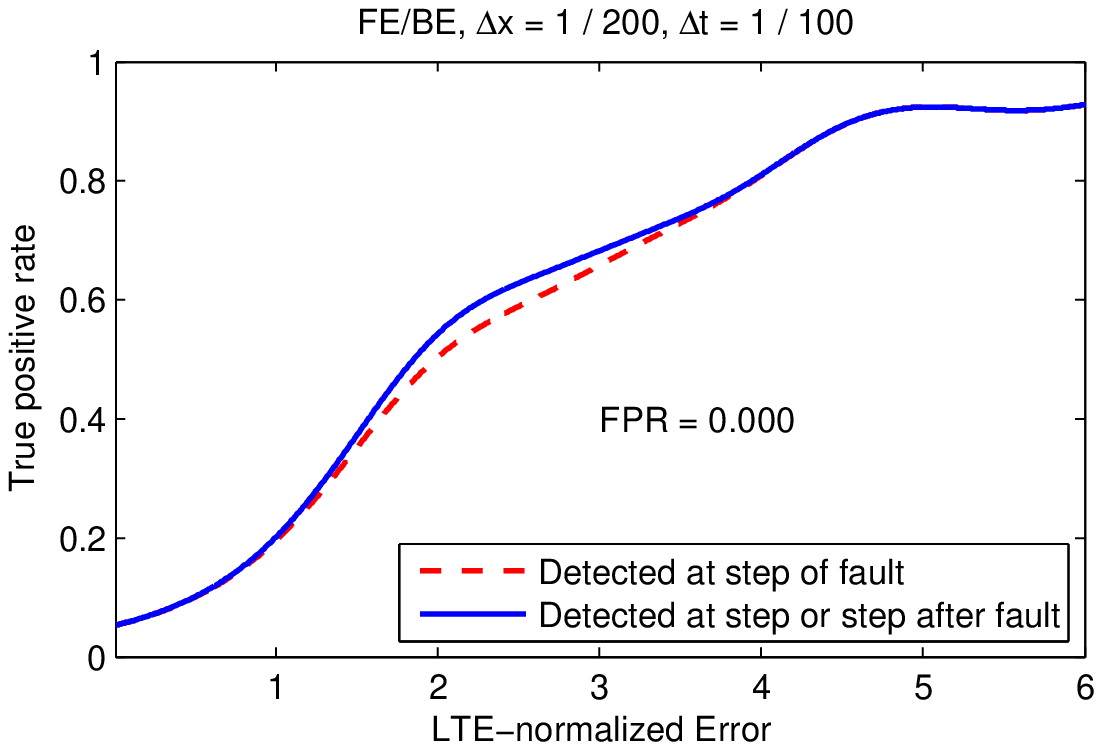}
  \includegraphics[scale=0.7]{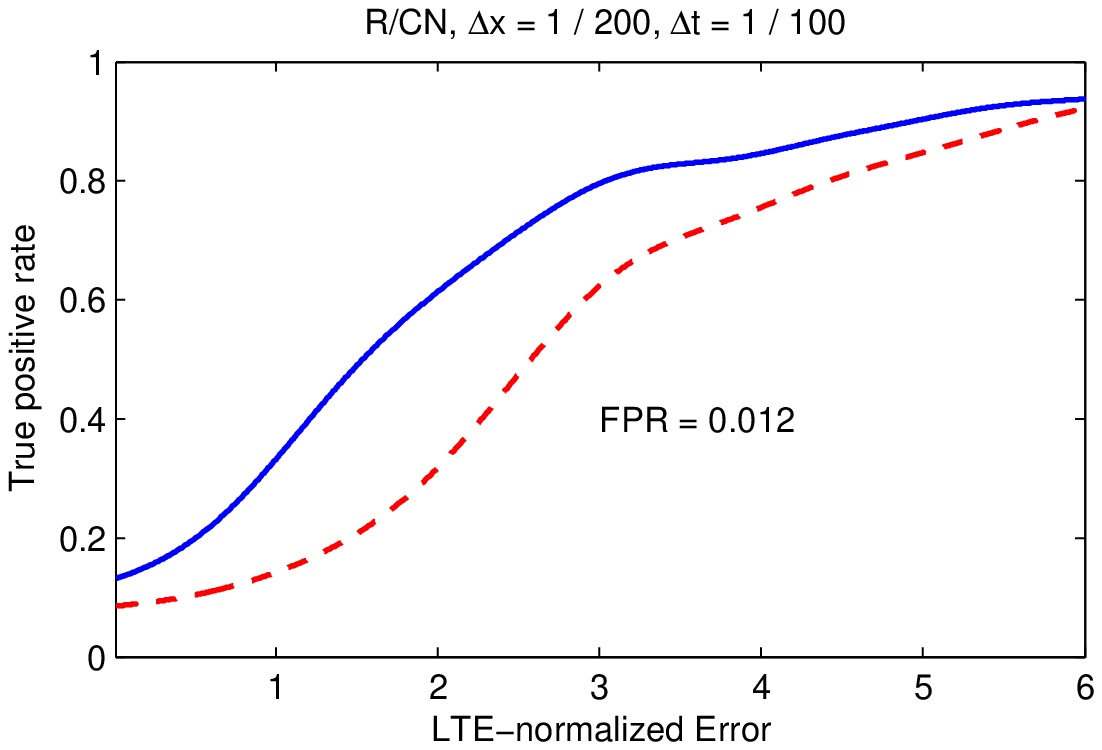} \\
  \includegraphics[scale=0.7]{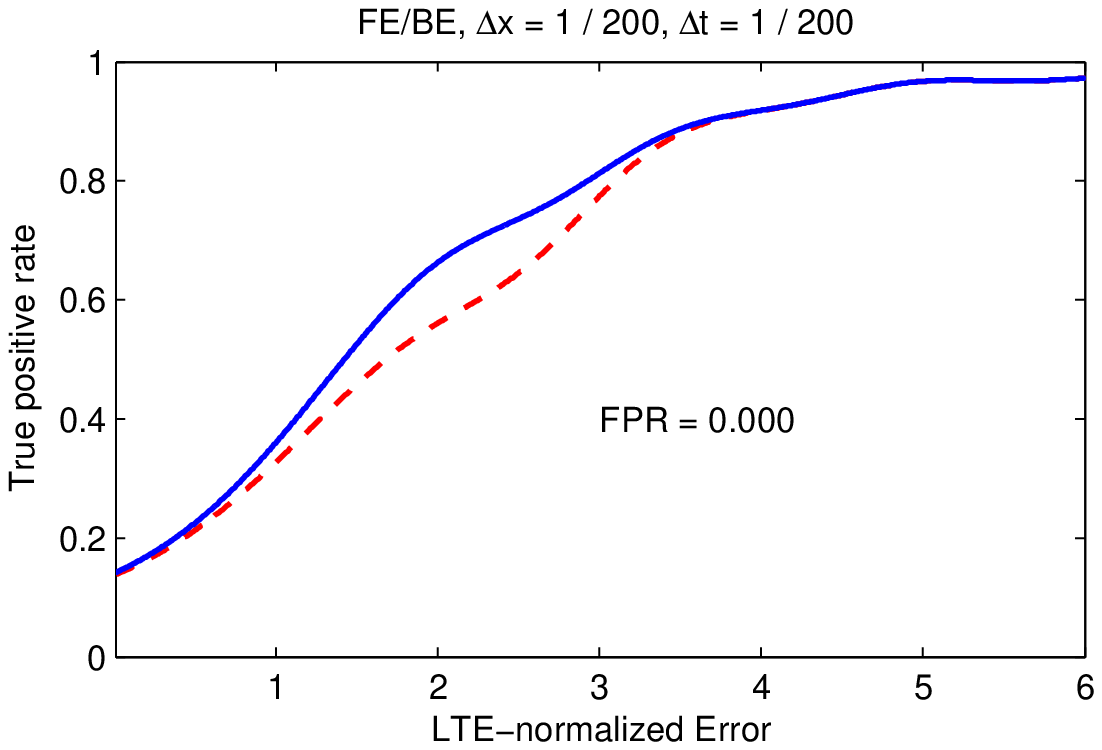}
  \includegraphics[scale=0.7]{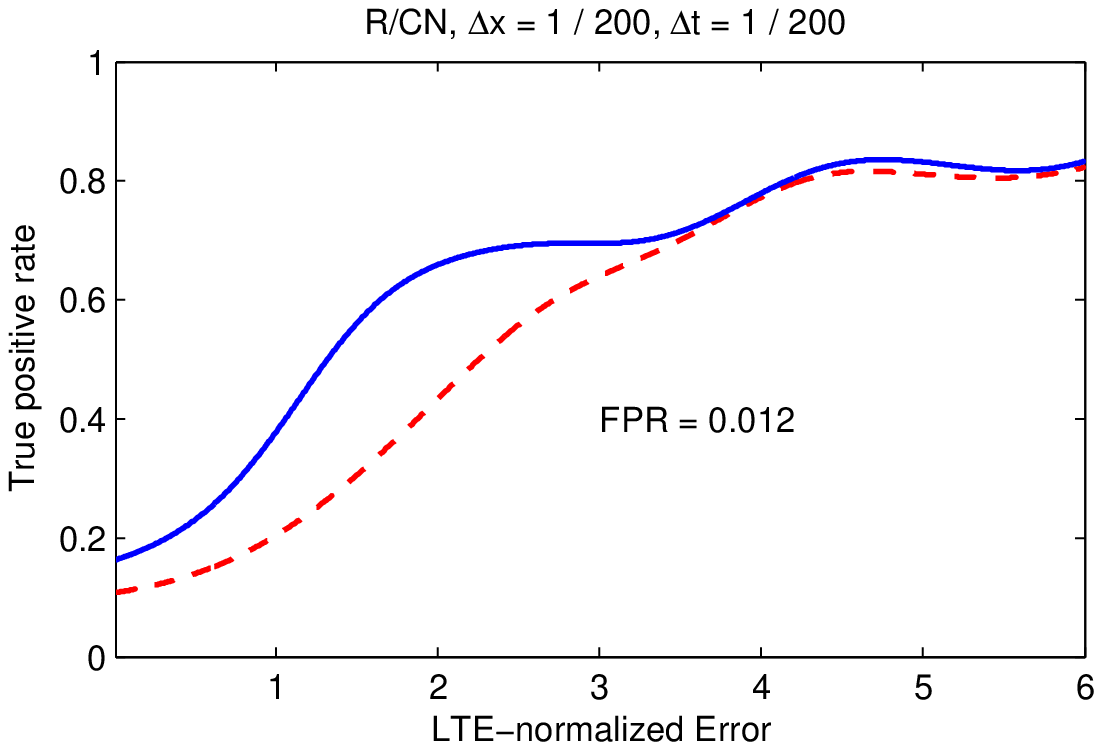} \\
  \includegraphics[scale=0.7]{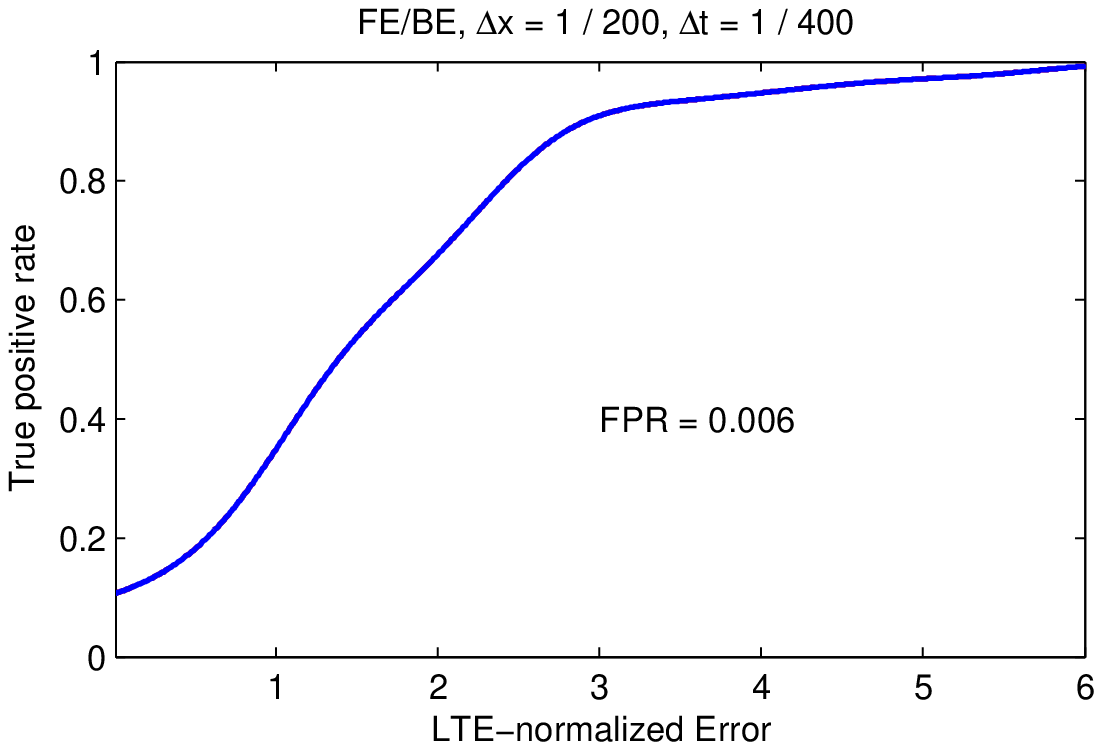}
  \includegraphics[scale=0.7]{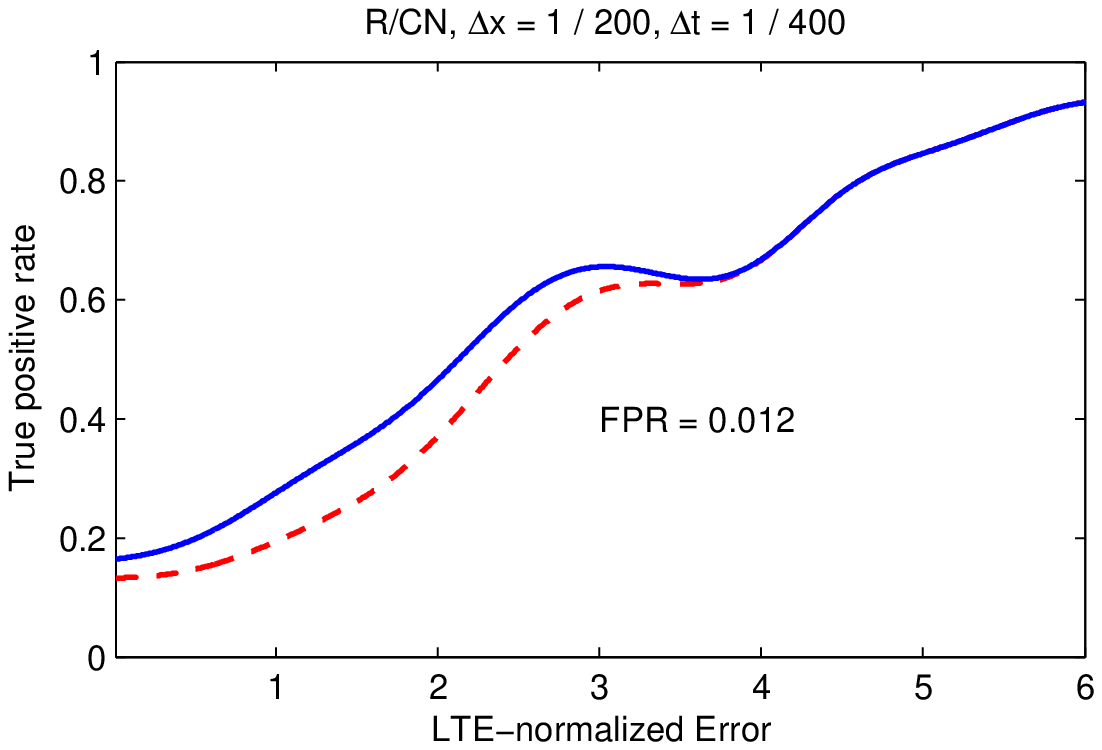}
  }}
  \caption{
  Detection results for Heat equation under Configuration 2 with $dt = 1 / 100$, $1 / 200$, and $1 / 400$.
  Faults are injected by multiplying a single random component of the right-hand-side of each linear equation
  solve by a normally distributed random variable with mean $1$ and variance 1e-6 (R/CN) or 5e-5 (FE/BE).
  }
  \label{fig:heat2}
\end{figure*}

\begin{figure*}[tbp]
  \centering
  \framebox{
  \parbox{470\unitlength}{
  \includegraphics[scale=0.7]{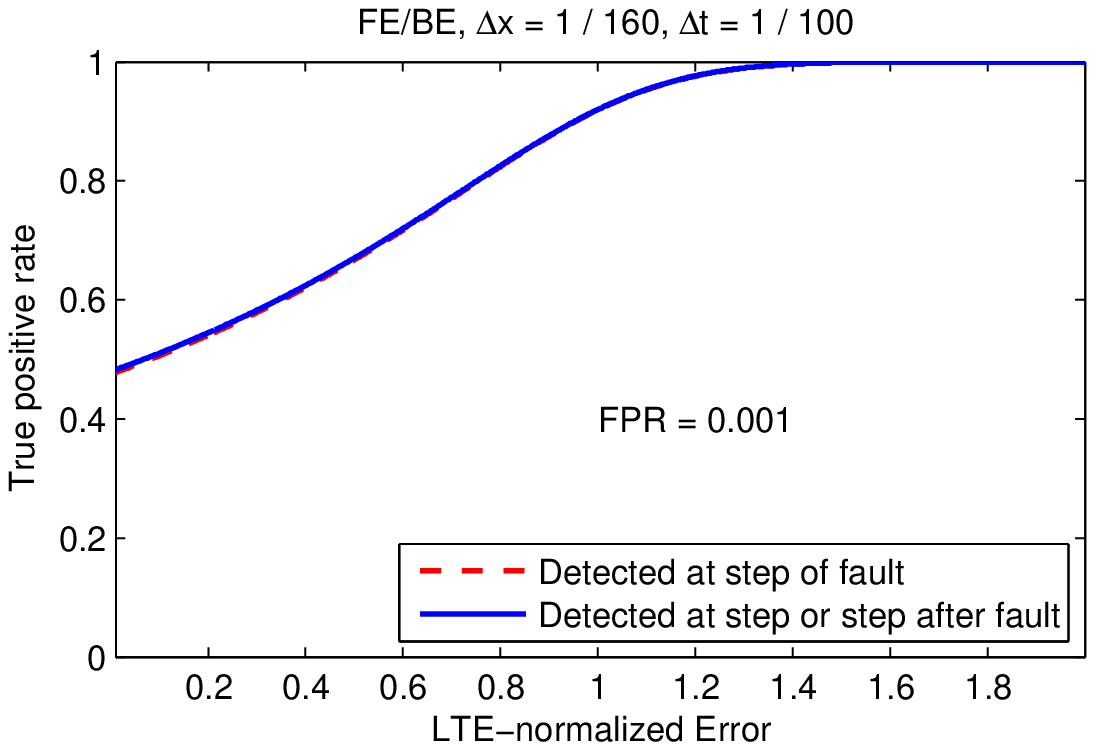}
  \includegraphics[scale=0.7]{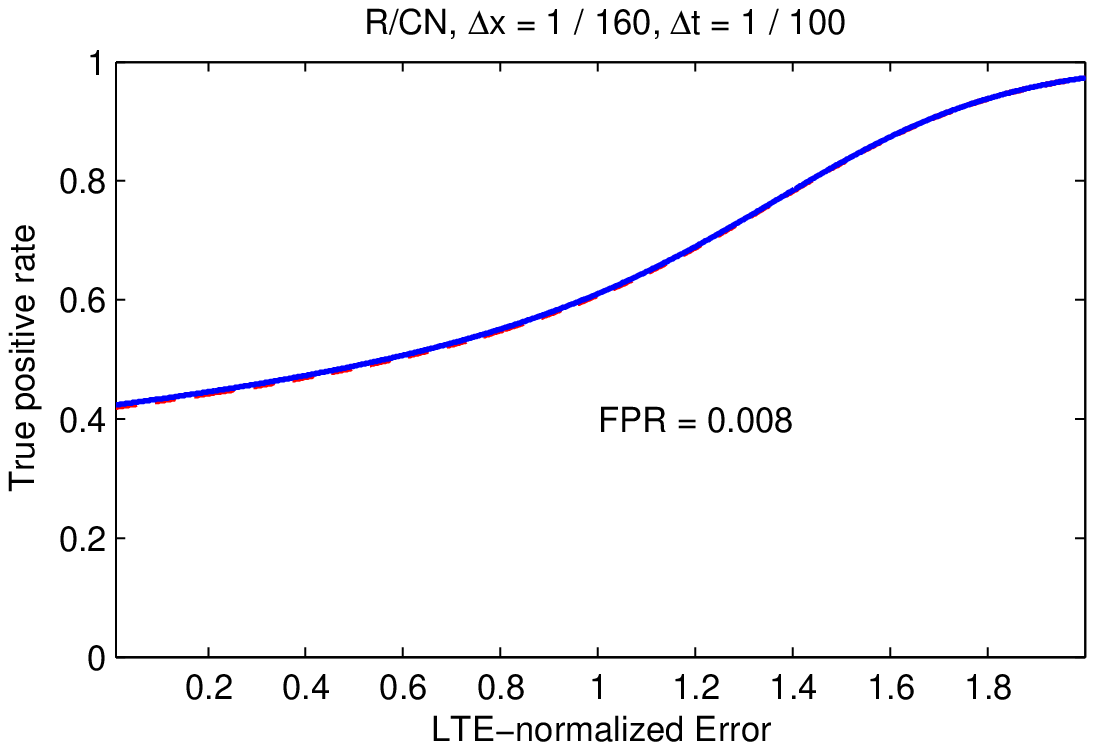} \\
  \includegraphics[scale=0.7]{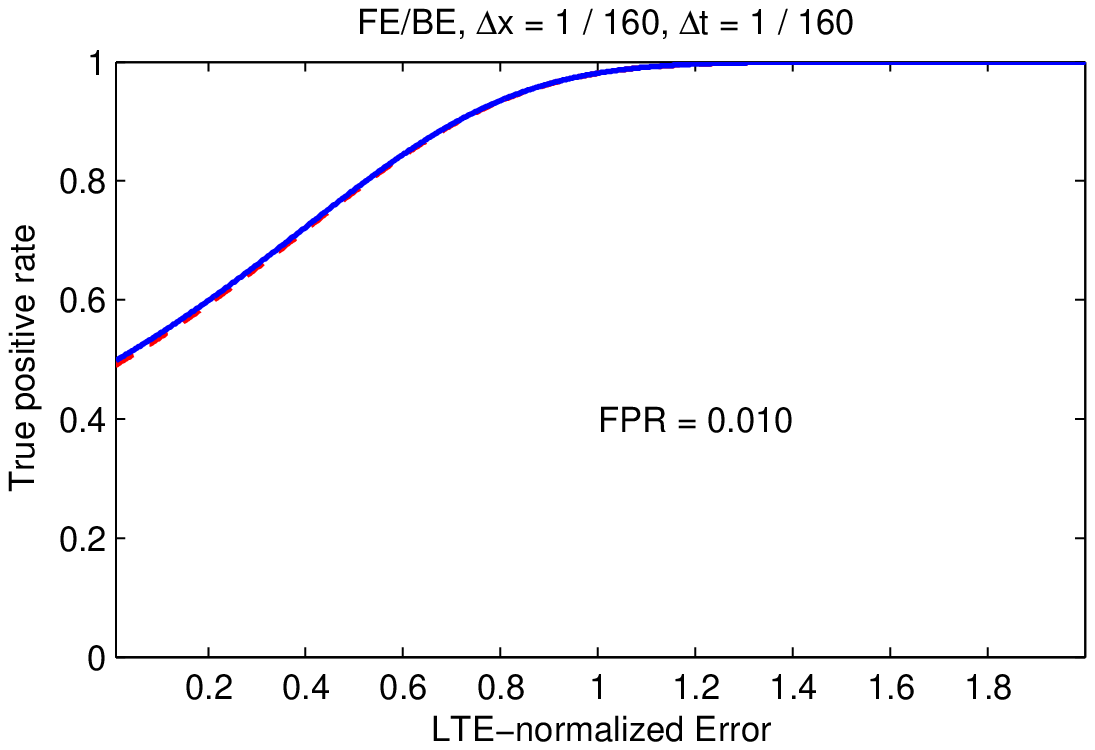}
  \includegraphics[scale=0.7]{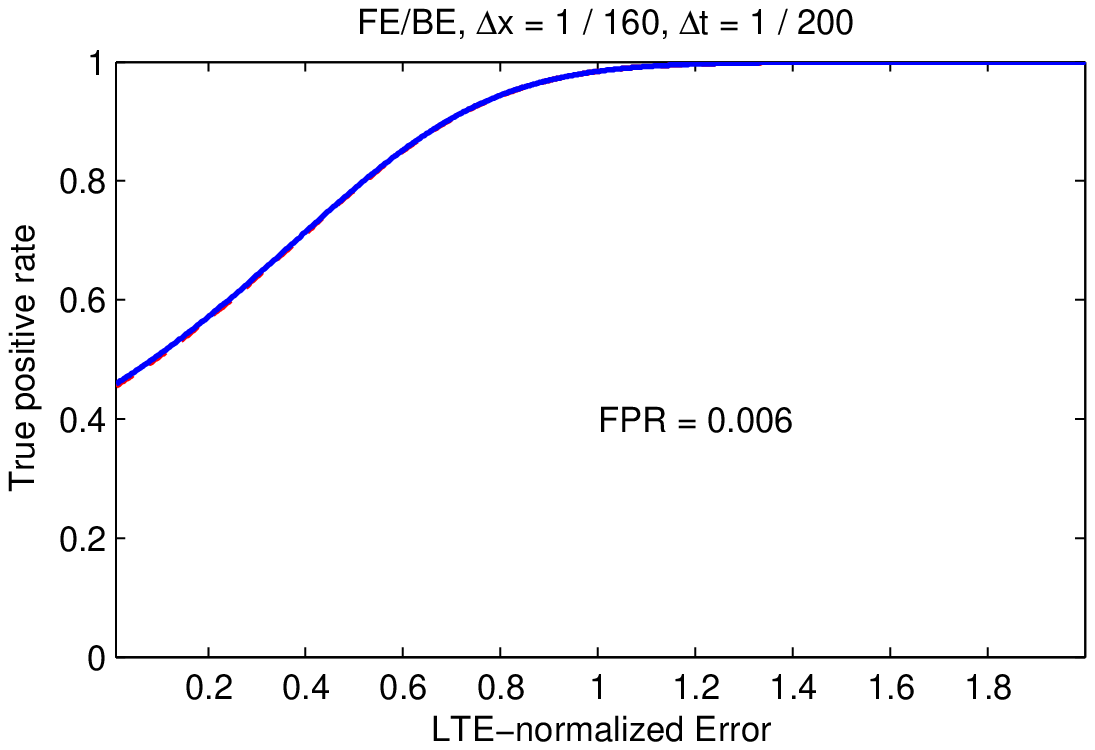} \\
  \includegraphics[scale=0.7]{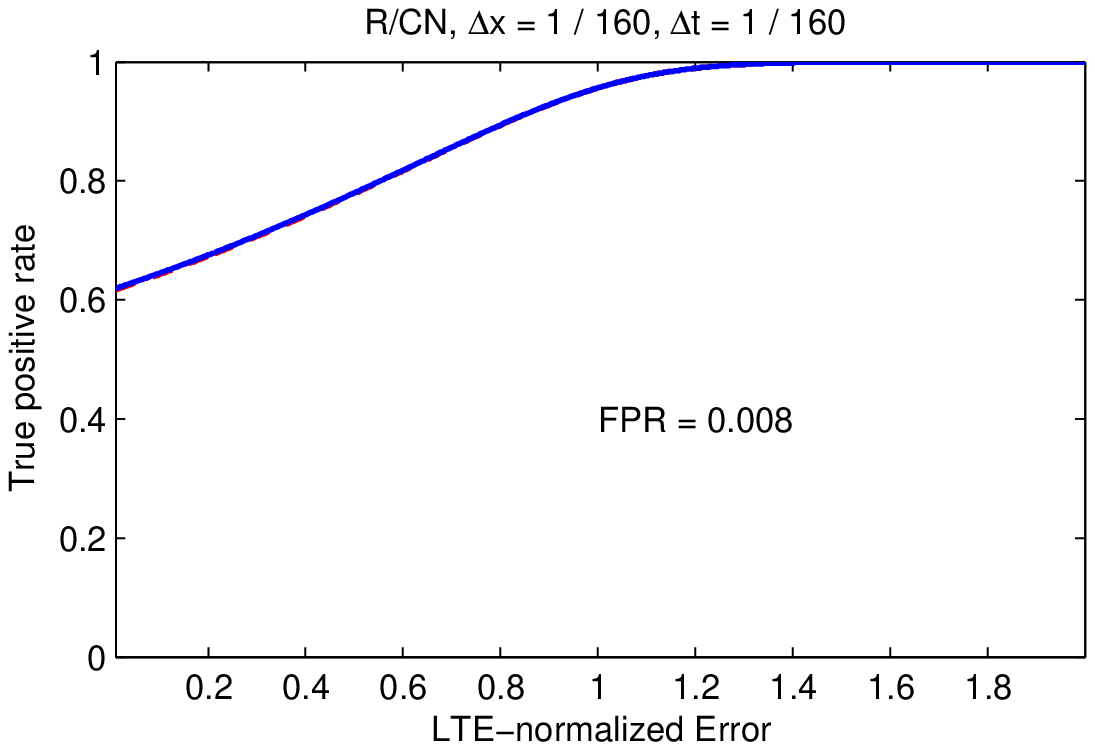}
  \includegraphics[scale=0.7]{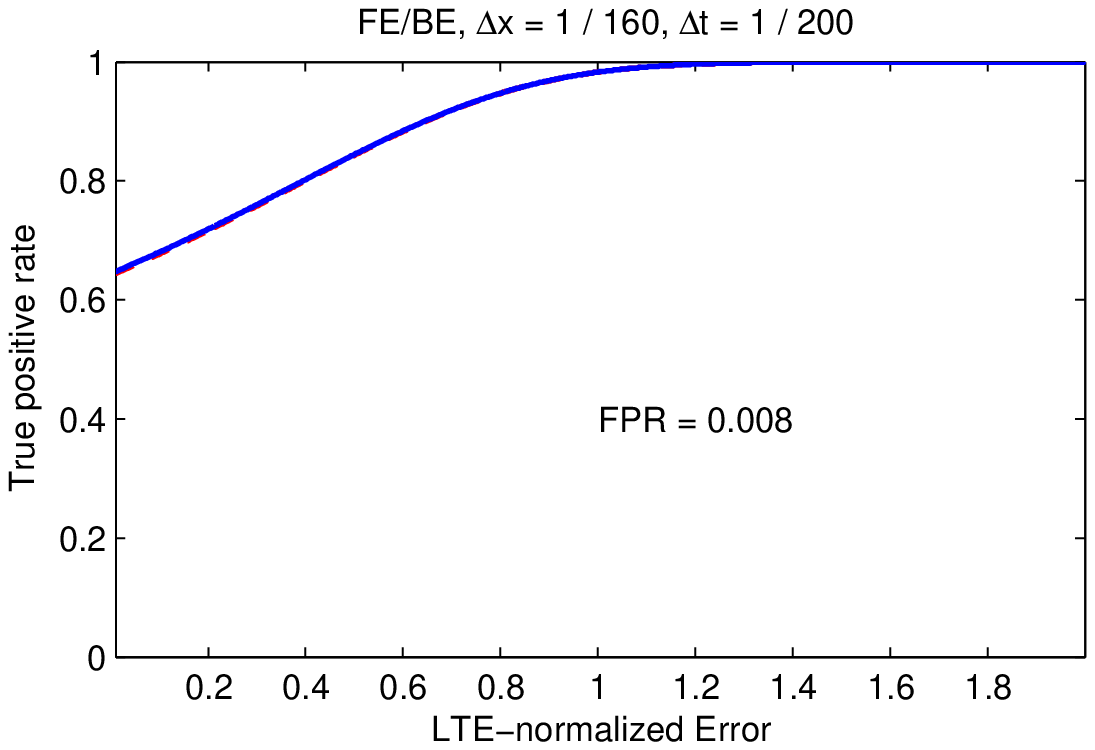}
  }}
  \caption{
  Detection results for Heat equation under Configuration 3 with $dt = 1 / 100$, $1 / 160$, and $1 / 200$.
  Faults are injected by multiplying a single component of the solution at the previous time step by a normally
  distributed random variable with mean $1$ and variance 1e-6 (Crank-Nicolson / Richardson) or 1e-4 (backward / forward Euler).
  }
  \label{fig:heat3}
\end{figure*}

We consider the heat equation (Equation~(\ref{eq:heat})) with homogeneous Dirichlet boundary conditions, for $x \in [0, 1]$, $t \in [0, T]$.
The $\A/\B$ formulations are the Richardson / Crank-Nicolson (R / CN) and
forward / backward Euler (FE/BE) schemes described in Section~\ref{sec:PDE_solvers}.
Table~\ref{tab:heat_configs} describes the three configurations of the heat equation used in our experiments.
For each configuration, we perform experiments with different time steps.

We consider a trial to be one call to the heat equation solver, which finds a numerical solution at
spatial points $0, \deltax, 2\deltax, \ldots, 1$ and temporal points $0, \deltat, 2\deltat, \ldots, T$.
We inject one fault per trial at a uniformly random step.
For Configuration 1, we corrupt a single component of the function evaluation of $q$ ($\sigma^2 = $1e-3 for R/CN and $\sigma^2 = $1e-1 for FE/BE).
For Configuration 2, we corrupt a single component in the right-hand-side of the implicit schemes' linear systems  ($\sigma^2 = $1e-6 for R/CN and $\sigma^2 = $5e-5 for FE/BE).
For Configuration 3, we corrupt a single component of the previous solution vector ($\sigma^2 = $1e-6 for R/CN and $\sigma^2 = $1e-4 for FE/BE).
The variances were chosen in order to generate errors with LTE-normalized error near one.
If the variances were much larger, nearly all errors would be detected and we would not see the relationship between TPR and LTE-normalized error.
Smaller variances were used for R/CN than FE/BE because the same type of corruption is more easily detected by R/CN.
This agrees with the case for the ODE solvers in Section~\ref{sec:vdp}, where higher-order $\A/\B$ schemes had better detection rates.

Figures~\ref{fig:heat1},~\ref{fig:heat2},~and~\ref{fig:heat3} plot TPR as a function of the LTE-normalized error.
The results illustrate several important features of the detection scheme.
First, errors with a large impact on the solution (large LTE-normalized error) are much more easily detected than
errors with a small impact (small LTE-normalized error).
Second, checking for an error one step after the fault occurs can significantly improve detection (see especially Figure~\ref{fig:heat1}).
In Section~\ref{sec:tardy}, we explore why this is true.
Third, the false positive rate (FPR) is small.
In many cases, no false positives are produced.
The largest FPR was only 1.2\%.
Fourth, we can detect several types of errors.
Finally, decreasing the time step either improves detection rates or keeps the detection rates the same.

%new
We note that the adaptive thresholding, described in Section~\ref{sec:error}, does not allow for a tradeoff between better TPR at the cost of a larger FPR or vice versa.
In a sense, the adaptive thresholding approximately finds the sweet spot where any anomaly that can be detected, is detected.
However, there are some anomalies that are so well ``disguised'', that they are indistinguishable from normal iterations, and only by allowing an extreme increase in FPR are we able to detect these.

\begin{comment}
We note that, in our experiments, we found that the adaptive thresholding (described in Section~\ref{sec:error})
makes it difficult to achieve a significantly better TPR at the cost of a larger FPR.
Increasing $\Gamma$ and decreasing $\gamma$ in Algorithm~\ref{alg:resilient} only slightly improves the TPR,
while the FPR can increase significantly if $\Gamma$ becomes too large or $\gamma$ too small.
\end{comment}

\subsection{Detection indicators for the heat equation}\label{sec:detect}

\begin{figure*}[tbp]
  \centering
  \fbox{
  \includegraphics[scale=0.6]{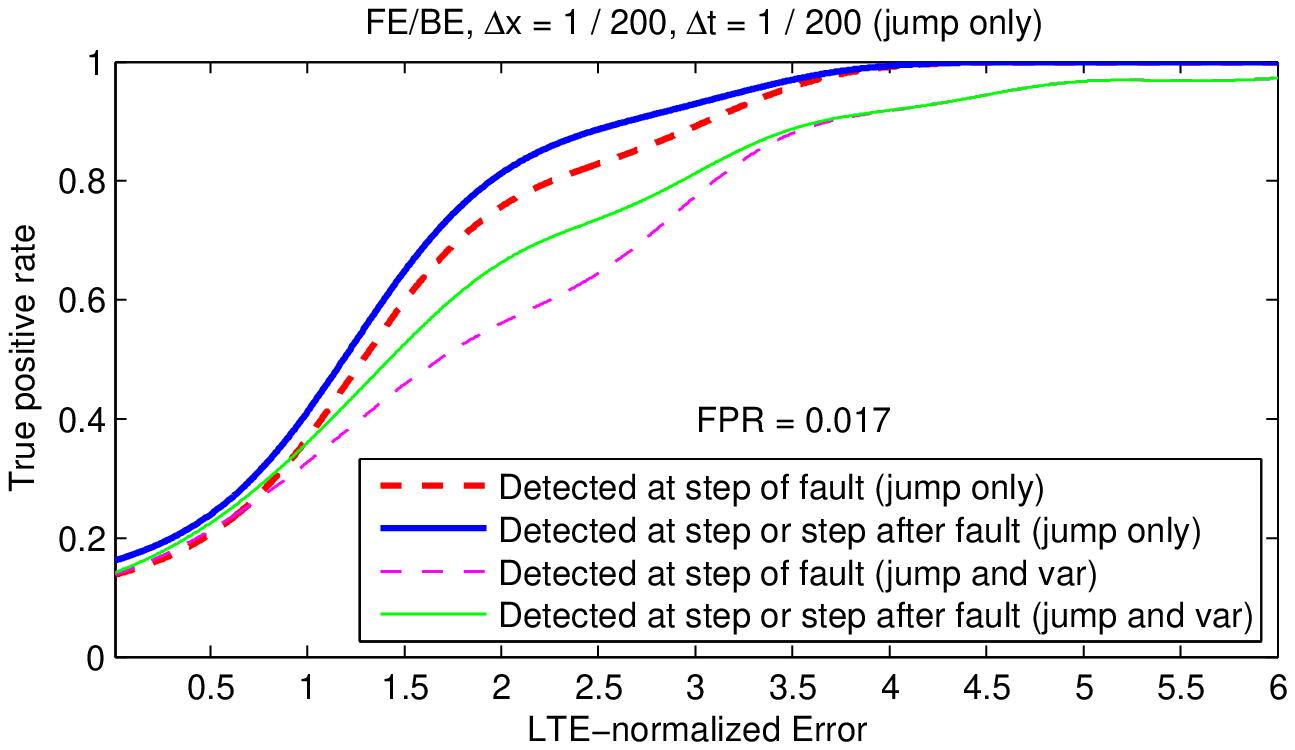}
  \includegraphics[scale=0.6]{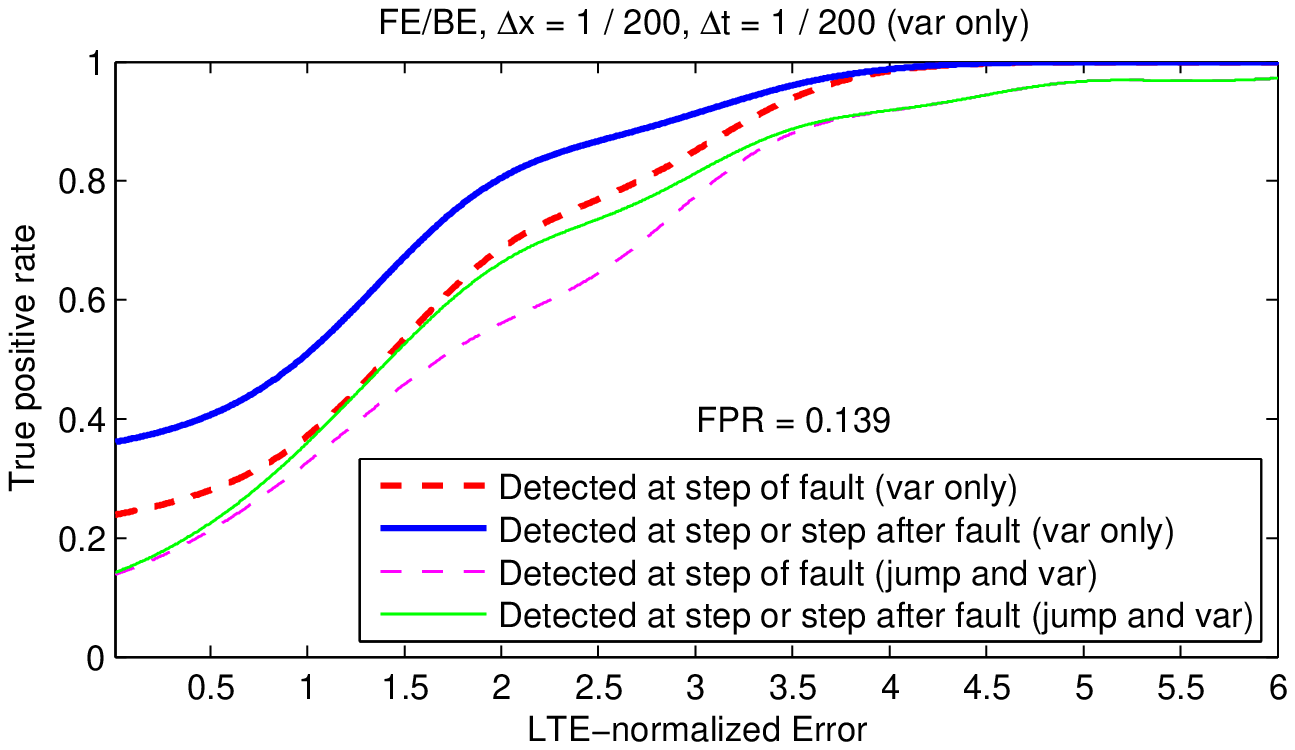}
  }
  \caption{
  Detection results for forward / backward Euler on the heat equation under Configuration 2 with $\deltat = 1 / 200$.
  The left and right plots compare using only the relative jump or only the variance change, respectively, to using both relative jump and variance change.
  The listed FPR is for using only the relative jump (left) or only the variance (right).
  When using both indicators, there are no false positives.
  Faults were injected by multiplying a single component of the solution at the previous time step by a normally distributed random variable with mean $1$ and variance 5e-5.
  Identical faults were injected for each type of detector.
  }
  \label{fig:heat_detect}
\end{figure*}

In Equations~(\ref{eq:relative_jump})~and~(\ref{eq:var_change}), we defined the indicators $J_n$ and $V_n$ used by the
detector's two-indicator strategy.
$J_n$ measured the jump in the differences in the sequence, and $V_n$ measured the change in variance of the differences.
We call these the \emph{relative jump} and the \emph{variance change}.

We now empirically explore the advantages of the two-indicator strategy over a single detection indicator.
Figure~\ref{fig:heat_detect} shows the detection results for FE/BE when using only the relative jump and
only the variance change under Configuration 2 of the heat equation with $\deltat = 1 / 200$.
We used the same injected faults as in Section~\ref{sec:numerical_heat}.
In other words, Figure~\ref{fig:heat_detect} shows the performance of the individual detectors when compared to the two-indicator strategy in Figure~\ref{fig:heat2}.

When the two-indicator strategy was used (Figure~\ref{fig:heat2}), no false positives were produced.
The FPR increases to 14\% when using only the relative jump and to 2\% when using only the variance change.
The TPR is only mildly increased when using only the relative jump and is mostly unchanged when using only the variance change.
This suggests that the FPR can be dramatically reduced by employing the two-indicator strategy.

\subsection{Tardy error detection with the heat equation}\label{sec:tardy}

\begin{figure*}[tbp]
  \centering
  \fbox{
  \includegraphics[scale=0.6]{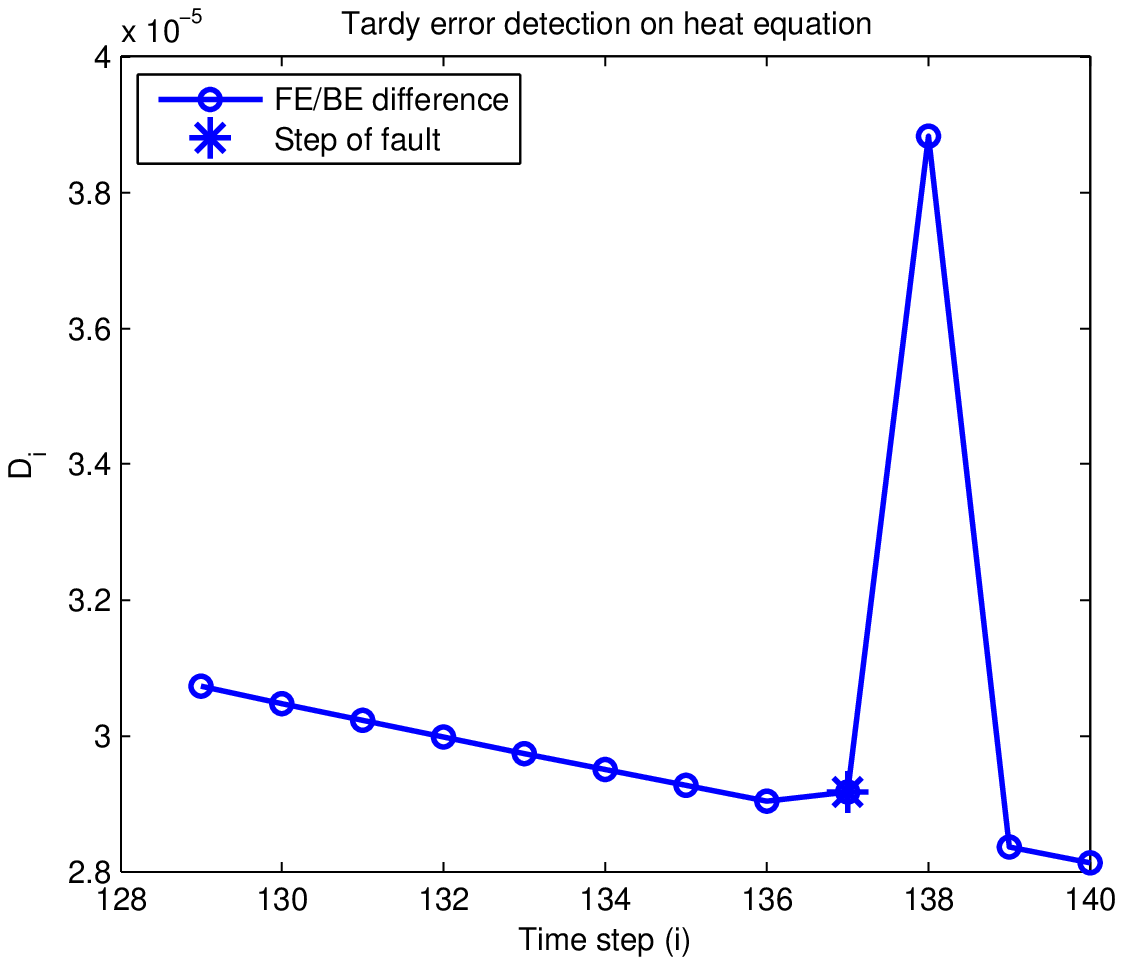}
  \includegraphics[scale=0.6]{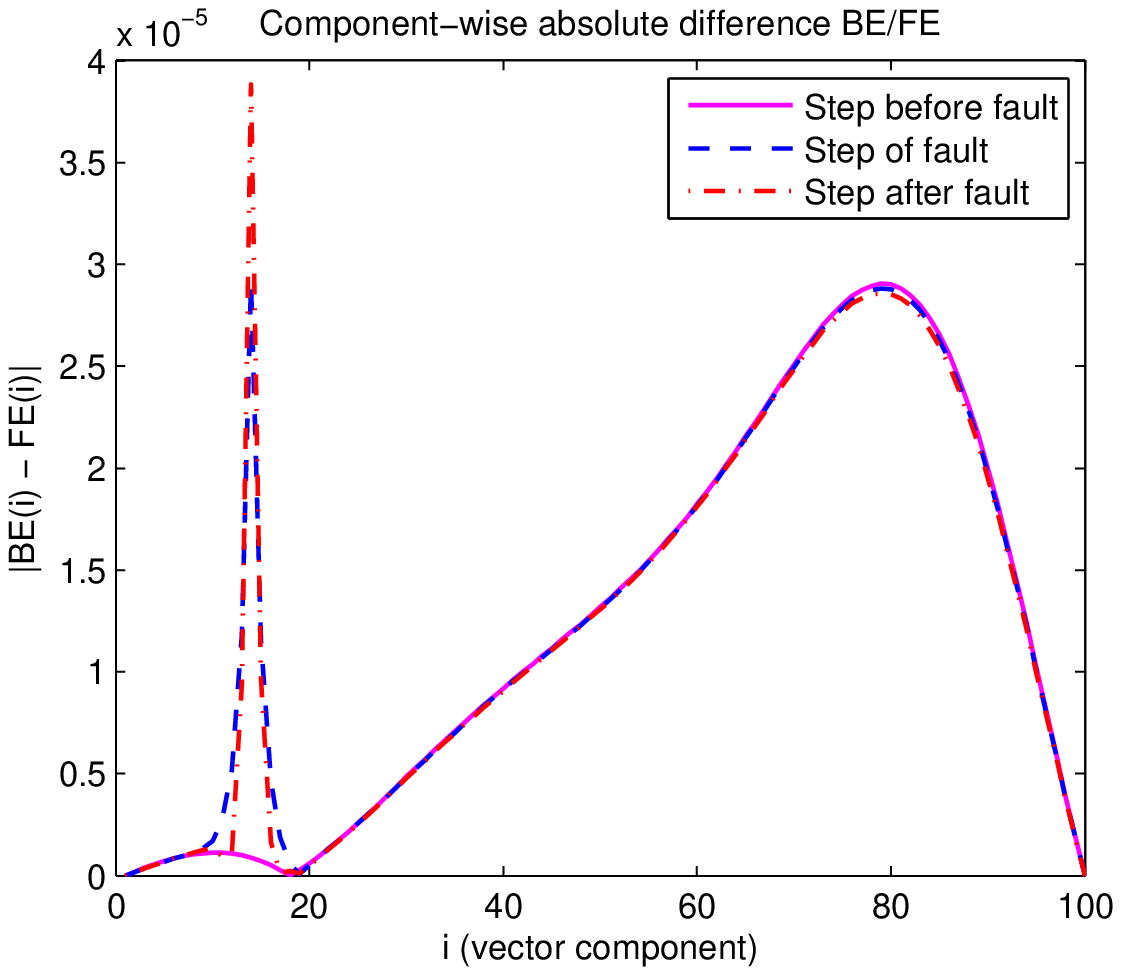}
  }
  \caption{
    The left plot shows the infinity-norm difference between forward and backward Euler schemes
    under Configuration 1 of the heat equation with $\deltat = 1 / 100$.
    At the step of the fault, there is a small jump in the difference, while at the step after the fault, there is a large spike.
    The right plot shows the component-wise absolute difference of forward and backward Euler solutions at the time steps before,
    of, and after the fault occurs.
    The fault corrupts the source term in the 13-th component, and a clear spike is seen at that location.
  }
  \label{fig:tardy_detection}
\end{figure*}

In some configurations of the heat equation, it is common to detect the error one time step \emph{after} the fault occurs (see, for example, Figure~\ref{fig:heat1}).
The left plot of Figure~\ref{fig:tardy_detection} shows $D_i$ for the FE/BE $\A/\B$ formulation near the point of the injected fault for one of the simulations in which tardy error detection occurs.
The right plot of Figure~\ref{fig:tardy_detection} shows the component-wise difference of the solution vectors near the time step of the fault (recall that the $D_i$ is the infinity-norm difference).
We see a small jump in $D_i$ at the step of the fault and a large spike the step after the fault.
The jump at the fault is not big enough for the detector to signal an error.
Some of the components of the difference between backward and forward Euler are naturally larger than others,
and the fault occurs at a spatial location where the differences tend to be small.

Why is the jump larger the step after the fault?
In this case, we are corrupting the source term $q$ in Equation~(\ref{eq:heat}).
At the fault, only backward Euler uses the corrupted evaluation of $q$ (forward Euler uses the value of $q$ from the \emph{previous} time step).
The implicit nature of backward Euler scheme forces the new solution to ``agree'' with the corrupted source term, and the result is a small jump in $D_i$.
At the step after the fault, forward Euler uses the corrupted source term.
Since forward Euler is explicit, the corrupted value is ``accepted'' and taken for a full time step.
This causes the large spike in the sequence of $D_i$ to occur at the step \emph{after} the fault.

The right plot of Figure~\ref{fig:tardy_detection} illustrates the advantages and disadvantages of using the infinity-norm.
On one hand, the errors tend to be localized spatially, and local spikes are easier to detect with the infinity-norm.
However, when the solution vector difference has different scales,
it is more difficult to detect faults that occur in spatial locations where the solution vector difference is smaller.
In general, we found that the infinity-norm worked better than the $1$-norm and $2$-norm.
We note that our general framework does restrict $D_i$ to consider only a single norm.
However, for simplicity, we chose a single norm for our experiments.
These results show that examing $D_i$ locally in space and in time can be beneficial, and this is an area of future work.

\subsection{Incompressible Navier-Stokes equations}
\label{sec:navierstokes}

The incompressible Navier-Stokes equations in two dimensions with no external forces are
\begin{align}
  u_t &= -(u^2)_x - (uv)_y + \frac{1}{Re}\LapOp u - p_x \nonumber\\
  v_t &= -(v^2)_y - (uv)_x + \frac{1}{Re}\LapOp v - p_y \nonumber.
\end{align}
Here, $u$ and $v$ are the velocity components, $Re$ is the Reynolds number, and $p$ is the pressure.

In our experiments, $\B$ is based on a simple projection method in \cite[Section~6.7]{Strang:2007} and open-source Matlab code \cite{Seibold:2008}.
The boundary is a square, and the boundary conditions are those of a driven cavity flow.

Let $U_{\B}^n$, $V_{\B}^n$, and $P_{\B}^n$ be the numerical solutions at the $n$-th time step and let $\deltat$ be the time step.
The overall structure of the update to $U_{\B}$ in an iteration of $\B$ is as follows:
\begin{enumerate}
\item Explicit (forward Euler like) handling of nonlinear terms:
\begin{equation}
\frac{U_{\B}^* - U_{\B}^n}{\deltat} = -((U_{\B}^n)^2)_x - (U_{\B}^nV_{\B}^n)_y, \nonumber
\end{equation}
where the subscripts denote centered difference.

\item Implicit solve for viscous term:
\begin{equation}
\frac{U_{\B}^{**} - U_{\B}^*}{\deltat} = \frac{1}{Re}\LapOp U_{\B}^{**} \label{eq:implvisc}
\end{equation}

\item Solve for the pressure correction:
\begin{equation}
\LapOp P^{n+1}_{\B} = \frac{1}{\deltat}\left(\left(U_{\B}^{**}\right)_x + \left(V_{\B}^{**}\right)_y\right) \label{eq:pressure}
\end{equation}

\item Update the solution:
\begin{equation}
U_{\B}^{n+1} = U_{\B}^{**} - \deltat (P_{\B}^{n+1})_x \label{eq:velupdate}
\end{equation}

\end{enumerate}
The update to $V_{\B}$ follows analogous steps.

For our experiments, we use extrapolation for the auxiliary scheme,
\begin{equation}
U^{n+1}_{\A} = 2U_{B}^{n} - U_{B}^{n-1}, \quad V^{n+1}_{\A} = 2V_{B}^{n} - V_{B}^{n-1} \nonumber
\end{equation}
along with the difference computation
\begin{equation}
D_{n} = \max(\max_{i, j} | (U^{n}_{\B})_{ij} - (U^{n}_{\A})_{ij} |, \max_{i, j} | (V^{n}_{\B})_{ij} - (V^{n}_{\A})_{ij} |) \nonumber.
\end{equation}
There is no restriction for the $\A/\B$ scheme to encompass the entire numerical method.
We could implement specialized $\A/\B$ schemes for each of the three steps in addition to or in place of the
extrapolation scheme.
An advantage of compartmentalized $\A/\B$ schemes is that we can detect an error early in the
iteration and avoid doing extra computation.
However, extrapolation is simple and demonstrates that detecting silent errors
in a nonlinear PDE system need not involve a lot of extra work.

Our experiments use the above projection method on the spatial domain $[0, 1] \times [0, 1]$ for $t \in [0, 2]$.
The discretizations are $\deltax = \deltay = 1/40$ and $\deltat = 1/100$.
We perform two simulations with different Reynolds numbers and different types of data corruption.
In the first simulation, $Re = 2000$, and we corrupt the previous solution, $U^n_{\B}$ ($\sigma^2 = $ 5e-1).
In the second simulation, $Re = 20$, and we corrupt the right-hand-side of the linear system in Equation~(\ref{eq:pressure}) ($\sigma^2 = 2$).
Each simulation consisted of 2,000 trials.
In each trial, the corruption occured at a single time step, chosen uniformly at random.
The entry in $U^n_{\B}$ and the entry in the right-hand-side of Equation~(\ref{eq:pressure}) were chosen uniformly at random.

\begin{figure*}
  \centering
  \fbox{
    \includegraphics[scale=0.6]{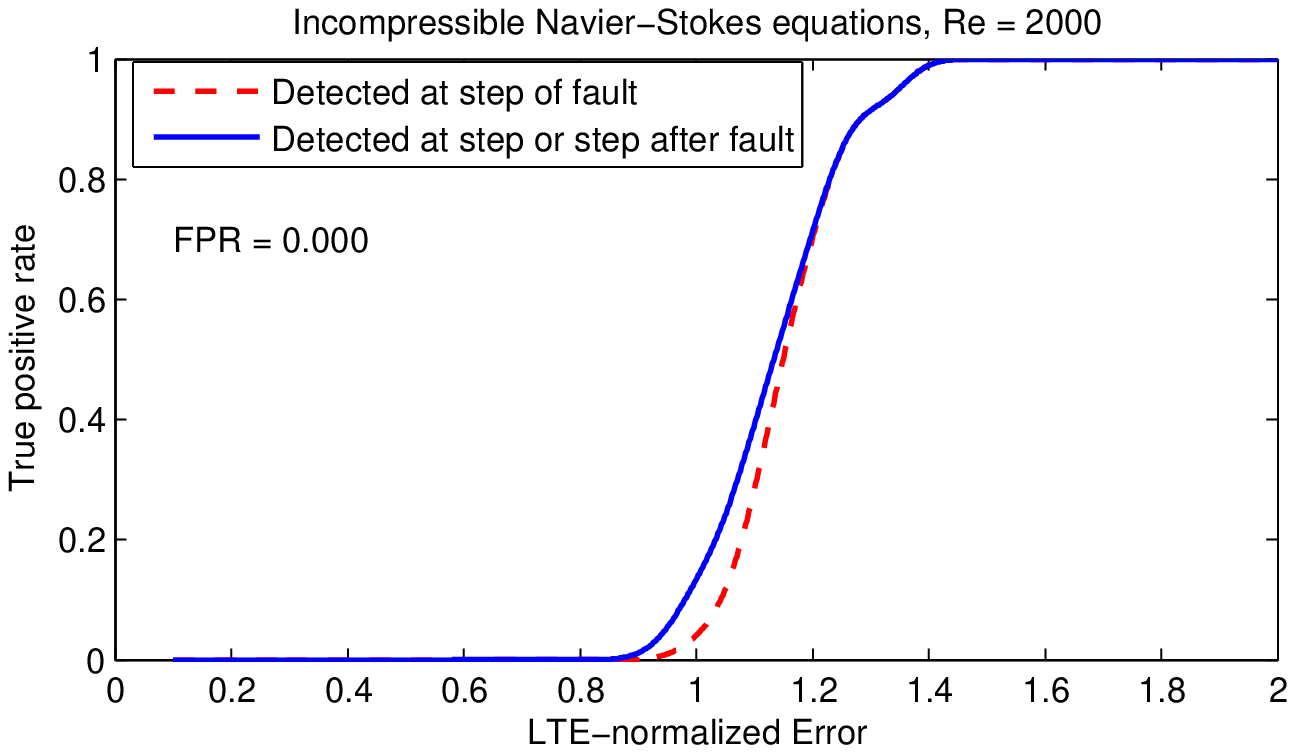}
    \includegraphics[scale=0.6]{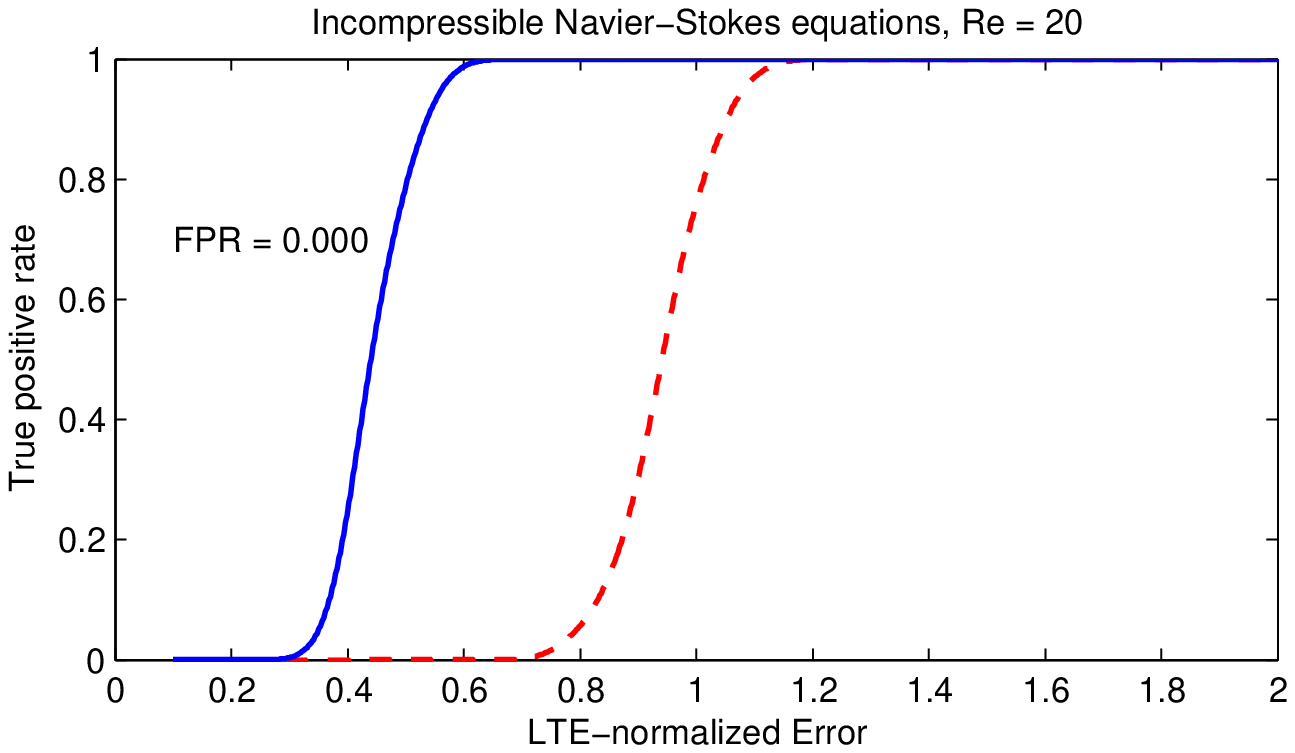}
  }
  \caption{
    Detector performance on incompressible Navier-Stokes equations with $Re = 2000$ (left) and $Re = 20$ (right).
    In the left plot, errors are introduced by multiplying a previous entry of the numerical solution of one velocity component ($U$)
    by a normal random variable with mean $1$ and variance 5e-1.
    In the right plot, errors are introduced by multiplying an entry in the right-hand-side of the linear system in Equation~(\ref{eq:pressure})
    by a normal random variable with mean $1$ and variance $2$.
  }
  \label{fig:ns_faults}
\end{figure*}

The TPR as a function of the LTE-normalized error is in Figure~\ref{fig:ns_faults}.
The results are similar to the behavior of the detector for the Van der Pol equation and the heat equation.
Errors with a large LTE-normalized error are easily detected, and the FPR is small.

  \section{Discussion}
\label{sec:conlusion}

By comparing the results of a base time-stepping scheme and an auxiliary scheme,
we are able to detect almost all significant errors.
The auxiliary scheme is readily available for standard ODE solvers such as Runge-Kutta and linear multistep methods,
as well as for PDE solvers for the Heat equation and the Navier-Stokes equations.
In simulations, our detection scheme successfully flags nearly all LTE-normalized errors above three,
while maintaining an overall false positive rate of less than 10\% (and in many cases, 0\%).
An important property of our detection scheme is that it is most successful detecting errors that have the largest impact on the solution.
We measure the impact by the LTE-normalized error.

One area of future work is a more formal analysis of the errors in the difference schemes and the sequence $D_i$.
We would like to say that a fault \emph{must} have occurred if some $D_i$ was above a computed threshold.
Typical error bounds are too loose with constants to be practical, so careful analysis is needed.

Further characterizations of silent errors is another area of future work.
First, it would be useful to detect \emph{what} caused the fault.
We can use data checksums to determine whether a previous solution was corrupted,
but determining if a function evaluation caused an error is more difficult.

Second, we would like to know \emph{where} the fault occurred.
For example, when perturbing an entry in the source term of the heat equation (the function $q$),
the heat is dissipated locally near the spatial point of perturbation.
The solution vector is then perturbed near (in space) to where the source term was perturbed.
We saw this phenomenon in Figure~\ref{fig:tardy_detection}, and it led to tardy error detection,
which we discussed in Section~\ref{sec:tardy}.
Thus, it is possible to detect where (physically) the fault occurred.
This is important for two reasons.
First, we can improve the performance of our error detector if we look for errors in space \emph{and} time.
Second, in parallel solvers, it is common for different spatial locations to be assigned to different processors.
By detecting the point in space where the fault occurred, we have an idea of which processor experienced a silent error.
In other physical simulations where perturbations cause local changes, we can apply the same idea.

The spatial location of the error is also potentially important in improving the sensitivity of the detector.   As shown in Figure~\ref{fig:tardy_detection}, it is possible that detectors that are local in both space and time may be a profitable extension of the approach taken here.

\section*{Acknowledgements}

We thank the US Department of Energy, which supported this work under Award Number DE - SC0005026.
Austin R. Benson is also supported by an Office of Technology Licensing Stanford Graduate Fellowship.
Sven Schmit is also supported by the Prins Bernhard Cultuurfonds.

%----------------------------------------------------------------------------------------
%	REFERENCE LIST
%----------------------------------------------------------------------------------------

  \bibliographystyle{apalike}
  \bibliography{main-bib}

%----------------------------------------------------------------------------------------

%\end{multicols}

\end{document}